\documentclass[amsmath,amssymb,aip,floatfix,reprint,cha]{revtex4-1}
\usepackage{amsfonts}
\usepackage{graphicx}% Include figure files
\usepackage{dcolumn}% Align table columns on decimal point
\usepackage{bm}% bold math
\usepackage[utf8]{inputenc}
\usepackage[T1]{fontenc}
\usepackage{mathptmx}
\usepackage{etoolbox}
\usepackage{siunitx}[v2]%=2021-04-09] % units
\DeclareSIUnit{\year}{a}
\usepackage{float} % control figure position?
\usepackage{multirow} % for tables
\usepackage[table]{xcolor}

\newcommand{\reco}{\text{R}_{\text{eco}}}
\newcommand{\R}{\mathbb{R}}

\renewcommand{\S}{Sec.} % refering to sections: changing from § to Sec. Change this as needed.

\makeatletter
\def\@email#1#2{%
 \endgroup
 \patchcmd{\titleblock@produce}
  {\frontmatter@RRAPformat}
  {\frontmatter@RRAPformat{\produce@RRAP{*#1\href{mailto:#2}{#2}}}\frontmatter@RRAPformat}
  {}{}
}%
\makeatother

\begin{document}

%\preprint{}
% \preprint{AIP/123-QED}

\title{%Orthogonal time-domain structure of dimensionality-reduction\\ in climatic time series analysis
%Time series analysis of vegetation dynamics data
% Extracting the seasonal cycle from vegetation dynamics data with dimension reduction methods
%Improved Time series Decomposition of Land-Atmosphere Fluxes with  Nonlinear Singular Spectrum Analysis
%Nonlinear Singular Spectrum Analysis detects seasonal cycle complexity in Land-Atmosphere Fluxes
%Nonlinear decomposition extracts harmonics in seasonal cycle of land-atmosphere fluxes
Nonlinear spectral analysis extracts harmonics from land-atmosphere fluxes
% or Nonlinear Time series Decomposition of Land-Atmosphere Fluxes 
} %Title of paper



\author{Leonard Schulz}
%\email[]%{leonard.schulz@posteo.net}
%leonard.schulz@ufz.de
%\thanks{}
\affiliation{Institute for Theoretical Physics, Leipzig University, Germany}
\affiliation{Institute for Earth System Science and Remote Sensing, Remote Sensing Centre for Earth System Research, Leipzig University, Germany}

\author{J\"urgen Vollmer}
\homepage{ORCID: \href{https://orcid.org/0000-0002-8135-1544}{0000-0002-8135-1544} }
\affiliation{Institute for Theoretical Physics, Leipzig University, Germany}

\author{Miguel D. Mahecha}
\homepage{ORCID: \href{https://orcid.org/0000-0003-3031-613X}{0000-0003-3031-613X} }
\affiliation{Institute for Earth System Science and Remote Sensing, Remote Sensing Centre for Earth System Research, Leipzig University, Germany}
\affiliation{German Centre for Integrative Biodiversity Research (iDiv) Halle-Jena-Leipzig, 04103 Leipzig, Germany}
\affiliation{Center for Scalable Data Analytics and Artificial Intelligence, ScaDS.AI Dresden/Leipzig, Germany}

\author{Karin Mora}
\email{karin.mora@uni-leipzig.de}
\homepage{ORCID: \href{https://orcid.org/0000-0002-3323-4490}{0000-0002-3323-4490} }
\affiliation{Institute for Earth System Science and Remote Sensing, Remote Sensing Centre for Earth System Research, Leipzig University, Germany}


\date{\today}

\begin{abstract} % max words = 250
% ABSTRACT: An article usually includes an abstract, a concise summary of the work covered at length in the main body of the article. It is used for secondary publications and for information retrieval purposes.
Understanding the dynamics of the land-atmosphere exchange of CO$_2$ is key to advance our predictive capacities of the coupled climate-carbon feedback system.
In essence, the net vegetation flux is the difference of the uptake of CO$_2$ via photosynthesis and the release of CO$_2$ via respiration, while the system is driven by periodic processes at different time-scales. 
%The net fluxes, obtained from the difference of uptake via photosynthesis and releases via respiration, are driven by periodic processes at different time-scales. 
The complexity of the underlying dynamics poses challenges to classical decomposition methods focused on maximizing data variance, such as singular spectrum analysis. 
Here, we explore whether nonlinear data-driven methods can better separate periodic patterns and their harmonics from noise and stochastic variability. 
We find that Nonlinear Laplacian Spectral Analysis (NLSA) outperforms the linear method and detects multiple relevant harmonics. 
However, these harmonics are not detected in the presence of substantial measurement irregularities. 
In summary, the NLSA approach can be used to both extract the seasonal cycle more accurately than linear methods, but likewise detect irregular signals resulting from irregular land-atmosphere interactions or measurement failures. Improving the detection capabilities of time-series decomposition is essential for improving land-atmosphere interactions models that should operate accurately on any time scale.
%However, this nonlinear method is not robust in the pressence of measurement irregularities. 
%Identifying the latter is, however, very important in carbon cylce research:
%unsuitable measurements must be identified such that they are not used e.g.~for model-calibration. 
%In summary, we recommend to adopt the NLSA approach to both extract the seasonal cycle and detect irregular signals due to irregular vegetation response or measurement failures.
\end{abstract}


\maketitle %\maketitle must follow title, authors, abstract and \pacs


% Lead paragraph:
% The first paragraph of the article should be the Lead Paragraph and contain the main points of the article, providing the “big picture” in a way that can be understood by non-specialist readers.
\begin{quotation}
    The analysis of oscillatory dynamics in time series data is a classical problem in celestial dynamics \cite{1998Gutzwiller,2000MurrayDermott,2013Laskar}, population dynamics \cite{2000Murray_Book,2017Barraquand-EtAl}, and chemical reactions \cite{1995Nicolis_Book} that still poses many research challenges.\cite{2001PikovskyRosenblumJuergenKurths_Book,2003KantzThomasSchreiber_Book,2018Ryashko,2018MaXuKurthsWangXu} 
    Understanding the land-atmosphere fluxes of CO$_2$ poses a complex challenge, because the CO$_2$ exchange is highly nonlinear and driven by different hydro-meteorological and biological factors.\cite{2010Hastings,bonan2015ecological,2017Barraquand-EtAl,2020WhiteHastings}
    These fluxes embody the multi-scale characteristics of hydro-meteorological dynamics, which are impacted instantaneously by variations in radiation e.g. due to varying cloud cover (spanning seconds to hours), daily alterations due to weather conditions, intra-seasonal shifts steered by synoptic variability, seasonal changes resulting from the Earth's revolution around the sun, and inter-annual differences triggered by intrinsic patterns of climate variability.
    % However, this understanding is crucial to enhance our Earth system models, thereby advancing our capacity to predict the future dynamics of the global carbon cycle \cite{hardouin2022uncertainty}.
    In the present paper, we analyze carbon flux measurements from the FLUXNET dataset \cite{baldocchi2001fluxnet} to test whether Nonlinear Laplacian Spectral Analysis (NLSA) outperforms linear decomposition methods when seeking to detect the seasonal cycle, and how it can be employed to warn about compromised data.
    Identifying the latter is crucial to avoid systematic errors in the calibration of Earth System models.\cite{mahecha2010comparing, hardouin2022uncertainty}
    %\todo{Watch out: I can not find Ref.~\citenum{richman2004sample}}
\end{quotation}

%########
%# Conventions
%# # section titles: Titles %start with a capital letter, but otherwise small.

%# General terms
%# # dimension reduction (NOT dimensionality reduction)
%# 

%#########################
%   S E C T I O N : Intro
%#########################
\section{Introduction}
\label{sec:intro}

\begin{figure}[b!]
    \centering
    \includegraphics[width=\linewidth]{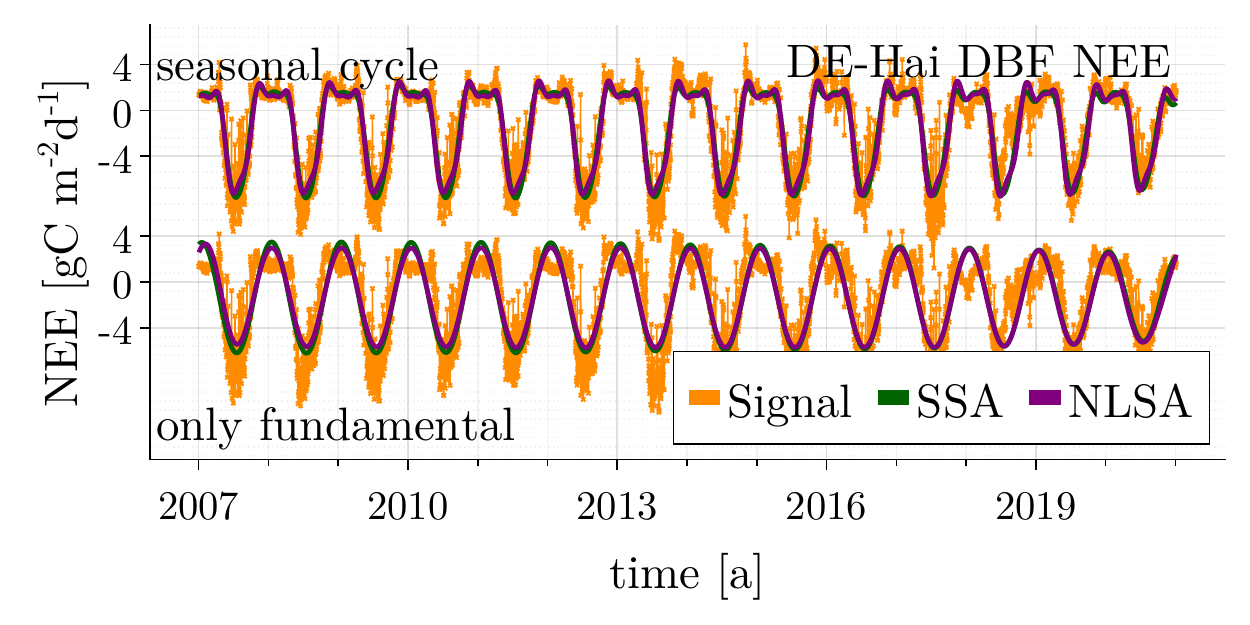}
    \caption{
    Net Ecosystem Exchange (NEE) measurements (signal, orange) of a deciduous broadleaf forest in central Germany (Hainich forest, DE-Hai, DBF). The seasonal cycle is approximated well by multiple harmonics (top row): four harmonics detected by NLSA, two harmonics by SSA. In contrast, the fundamental oscillation (lower row), typically detected by linear spectral analysis fails to accurately represent the variability during the summer period. 
        %In this example, both methods perform well, i.e. the NLSA and SSA curves are very similar.
        }
    \label{fig:motivation}
\end{figure}

% basics of carbon cycle
The continuous exchange of carbon dioxide (CO$_2$) between terrestrial ecosystems and the atmosphere is one of the key regulators for atmospheric CO$_2$ concentrations. In the presence of active vegetation, CO$_2$ is absorbed over land during during the daytime via photosynthesis. At the ecosystem and global scale this process determines the largest CO$_2$ flux from the atmosphere to the land and is referred to as gross primary production (GPP). In parallel, CO$_2$ is constantly released back to the atmosphere through various respiration activities, collectively known as ecosystem respiration ($\reco$).
%(R$_{\text{eco}}$).
The difference between these two counteracting fluxes defines the net ecosystem exchange (NEE). \cite{schulze2006biological, luyssaert2009toward} 
When integrated over time, NEE determines whether an ecosystem acts as a natural carbon dioxide sink (NEE$<0$) or source (NEE$>0$). 
%\todo{Here we introduce NEP. Later we always consider NEE. Better drop acronym NEP, and introduce NEE here.}
An example for its temporal evolution is shown in Fig.~\ref{fig:motivation}.
Given the rapid increase of human-induced CO$_2$ emissions into the atmosphere in our current era, and the uncertain fate of the terrestrial carbon sink, it is paramount to understand how these natural fluxes vary across different locations and over time 
to enhance our Earth system models, thereby advancing our capacity to predict the future dynamics of the global carbon cycle.\cite{mahecha2010comparing, hardouin2022uncertainty} 

% Understanding the land-atmosphere fluxes of CO$_2$ poses a complex challenge, because the processes underlying GPP and $\reco$ are highly nonlinear and driven by different hydro-meteorological and biological factors \cite{bonan2015ecological}. These fluxes embody the multi-scale characteristics of hydro-meteorological dynamics, which are impacted instantaneous variations in radiation e.g. due to varying cloud cover (spanning seconds to hours), daily alterations due to weather conditions, intra-seasonal shifts steered by synoptic variability, seasonal changes resulting from the Earth's revolution around the sun, and inter-annual differences triggered by intrinsic patterns of climate variability. Moreover, ecosystems exhibit internal resilience to hydro-meteorological changes through their buffering mechanisms. For instance, they can maintain water for an extended duration, contributing to heat storage, and exhibit vegetation responses that interact with the environment, thereby introducing additional modes of variation \cite{baldocchi2020eddy}.

% methodological approaches: our previous work
The multi-scale nature of land-atmosphere fluxes can be modeled analytically with a stochastic model\cite{pappas2017ecosystem} and has also been empirically corroborated using non-destructive, prolonged measurements of NEE, GPP, and $\reco$ that span across multiple years without interruption.\cite{mahecha_characterizing_2007, vargas2010multiscale, stoy2009biosphere} Studies aiming to extract modes of climate variability on multiple time scales use different decomposition methods, such as Fourier decomposition \cite{linscheid2020towards}, Wavelet spectra \cite{stoy2009biosphere}, Singular Spectrum Analysis\cite{mahecha_characterizing_2007, biriukova2021performance}, or empirical mode decomposition. \cite{mahecha_identifying_2010} While these methods are efficient, they neglect the geometry of the data manifold. Arguably, this limits the capability to capture features with low variance but of dynamical importance. \cite{giannakis_time_2011,giannakis2012nonlinear,giannakis2013nonlinear,berry2013timescale} 
%their basis are orthogonal functions which could potentially limit their capability to fully capture the intricacy of the underlying signals. 
% Being confined to linear basis functions could result in imprecise time-scale separations, as the processes under investigation are not precisely isolated and spectral leakage might mask the patterns of interest. However, considering that such decomposition methods may form the basis for estimating process variables on specific time scales  \cite{mahecha_identifying_2010, liu2020soil, wu2021diagnosing}, precision of the time-scale seperation is of utmost importance. Biogeochemical models parameterized based on such estimates might be integrated in Earth system models for future scenarios, making accuracy even more critical.
Moreover, spectral leakage might mask the patterns of interest as the processes under investigation are not precisely isolated. 

% Approach in this study
In the present study, we analyze time series of vegetation fluxes using two dimension reduction methods for temporal pattern extraction: the linear method Singular Spectrum Analysis\cite{broomhead1986extracting} (SSA) and the nonlinear method Nonlinear Laplacian Spectral Analysis\cite{giannakis2012nonlinear} (NLSA). Our goal is to improve the understanding of the underlying dynamics, specifically harmonic components. We investigate whether the nonlinear method can substantially improve the extraction of periodic patterns in land-atmosphere fluxes of CO$_2$. %(construction of the seasonal cycle in the upper panel of Fig.~\ref{fig:motivation}, rather than only the fundamental mode in the lower panel).
% in order to improve our understanding of the underlying dynamics. Specifically, we compare a novel %decomposition technique, which incorporates nonlinear dynamics, to the classical SSA approach. 
%Our ambit is to better capture the intricate interplay of hydro-meteorological and biological factors that drive the fluxes of CO$_2$, water, and energy between terrestrial ecosystems and the atmosphere. 
% The most crucial step is extracting the seasonal component from the remaining modes of variability, whereby the former is identified by its fundamental oscillation of period one year. 
% We expect that t
% The representation of the non-trivial seasonal dynamics can be substantially improved by NLSA.
% including harmonic oscillations. 
% Especially the NEE time series is periodic with more complex dynamics during the summer period due to the contrasting responses of the contained photosynthetic and respiratory fluxes as one can see in Fig.~\ref{fig:motivation}. 
By periodic patterns we mean higher order harmonic oscillations \cite{Bacastow1985seasonal} from which we construct the seasonal cycle. The fundamental oscillation does not capture the intricate interplay of hydro-meteorological and biological factors that drive the fluxes of CO$_2$, water, and energy between terrestrial ecosystems and the atmosphere at specific locations, Fig.~\ref{fig:motivation}.
%, which illustrates complex periodic dynamics during the summer period due to the contrasting responses of the contained photosynthetic and respiratory fluxes. NLSA allows us to capture this by representing the seasonal cycle with several harmonics.  
%into account,
% which yields a better approximation. 
% To detect these additional harmonic oscillation nonlinear methods are better suited. 

The paper is structured as follows: In \S~\ref{sec:MatAndMeth} we introduce the data, including the basic data measurement principle, the data selection, and preprocessing steps. We also summarize our methodological approach and the two data-driven methods used in this study: the linear method SSA and the nonlinear method NLSA. 
The results are presented in \S~\ref{sec:results}, where we highlight the challenges posed by irregular time series and provide an overview of the methods' performance across measurement variables and measurement sites.
In \S~\ref{sec:discuss}, we first discuss the implications of constructing the seasonal cycle from the harmonics obtained with SSA and NLSA.
Next, we compare the methods' performance in detecting harmonics. Finally, we outline the potential of our approach to study biosphere dynamics, meteorological dynamics, and their interactions. Section~\ref{sec:conlusion} concludes our study.

\section{Materials and Methods}
\label{sec:MatAndMeth}    
    \subsection{Data: Land-Atmosphere Fluxes}
    \label{sec:data}

    \begin{figure}[htb!]
        \centering
        \includegraphics[width=0.84\linewidth]{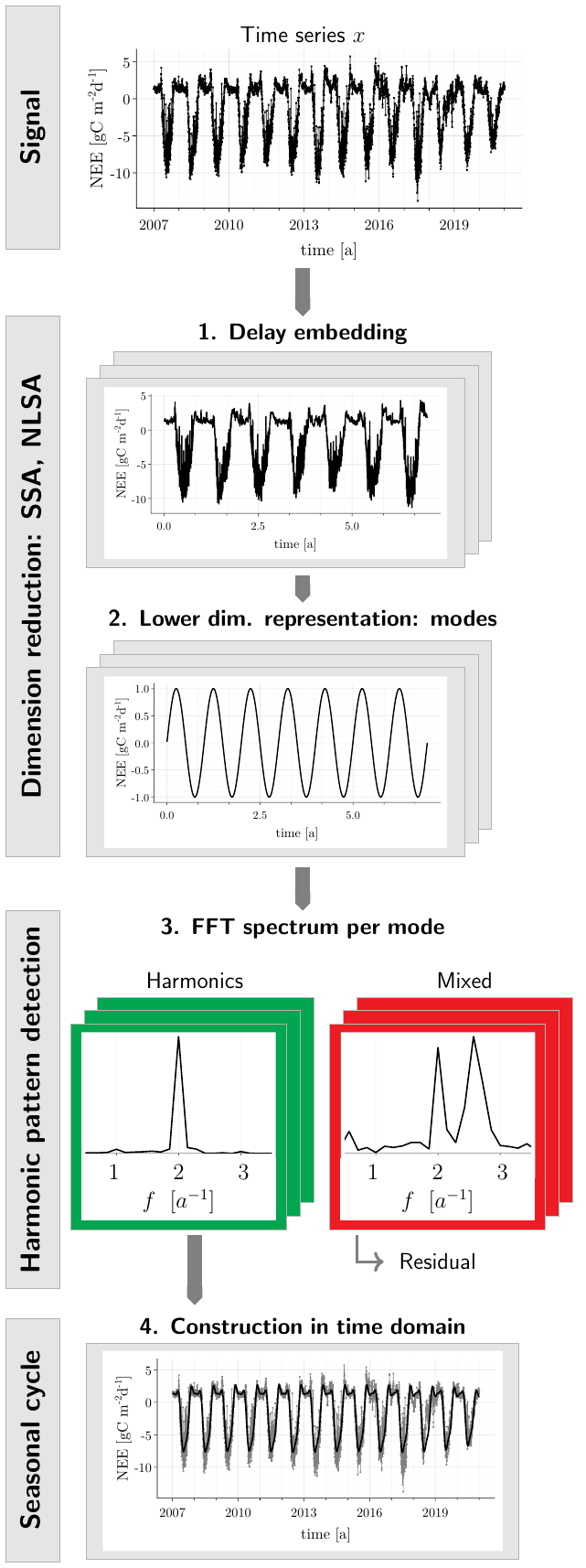}
        \caption{
        Workflow.
        The time series $x$ (Signal) is embedded into a higher-dimensional space using Takens' delay coordinates, which yields the data matrix $X$. The dynamics are approximated by modes computed with dimension reduction techniques SSA (linear) or NLSA (nonlinear). These modes are further analyzed using FFT to identify modes representing pure harmonic oscillations. Finally, the identified harmonic oscillations are used to construct the seasonal cycle.
        }
        \label{fig:workflow}
    \end{figure}
    
    Land-atmosphere fluxes, such as CO$_2$, are measured using the eddy covariance technique, \cite{aubinet1999estimates,baldocchi2001fluxnet,baldocchi2020eddy} a widely used method in atmospheric science, ecology, and environmental monitoring. This technique is primarily employed to measure the exchange of gases, such as carbon dioxide and water vapor, between ecosystems (e.g., forests, grasslands) and the atmosphere.
    The eddy covariance method involves measuring the fluctuations in the vertical wind speed and the concentration of gases at high frequencies (typically $10$–$\qty{20}{\hertz}$). These measurements are taken on a flux tower over an ecosystem of interest. The tower is equipped with sensors to capture the 3D variations in wind speed and gas concentrations, enabling to extract the turbulent exchange fluxes between the land-surface and atmosphere.
    These continuous, high-frequency measurements are available at various ecological in-situ monitoring networks, for instance via ICOS (\url{www.icos-cp.eu}) in Europe or FLUXNET at the global scale (\url{fluxnet.org}). 
    
    In this study we analyze time series from nine sites on the European continent, Table~\ref{table:sites}, from the ICOS 2020 Warm Winter dataset.\cite{WarmWinter2020Team2022} Per site 
    we analyze seven variables individually:
    GPP $(\unit{g C m^{-2} d^{-1}})$, 
    $\reco$ night $(\unit{g C m^{-2} d^{-1}})$, 
    NEE $(\unit{g C m^{-2} d^{-1}})$, as well as 
    air temperature (TA) $(\unit{\degreeCelsius})$, 
    short wave radiation (SW) $(\unit{W m^{-2}})$, 
    soil temperature (TS) $(\unit{\degreeCelsius})$, and 
    soil water content (SWC) (\%). 
    We focus on forest sites, as these are affected the least by short term human interventions and use data from two deciduous broadleaf forests (DBF), five evergreen needleleaf forest (ENF), and two mixed forest (MF). We seek the longest temporal measurements to best capture longer cyclic dynamics. These criteria are satisfied at a total of nine eddy covariance sites, where data was collected throughout the period 01.01.2007--31.12.2020 with a sampling rate of one measurement per day. 
    %located in Belgium (BE), Denmark (DK),  Germany (DE), Finland (FI), Italy (IT), Russia (RU), and Switzerland (CH).
    %moved to methods section: We also show quality flags (QF) at discrete time points. They indicate the confidence of measurement quality or in the gap filling process in \%, i.e. a gray scale, where no confidence is indicated by black. 
    
\begin{table}[h!]
    \caption{Description of the nine selected sites on the European continent, which are sorted by land cover: evergreen needleleaf forest (ENF), mixed forest (MF), deciduous broadleaf forests (DBF).
    }
    \centering
    \begin{ruledtabular}
    \begin{tabular}{c c c c}
    Site ID	&	Country	&	Site Name	&	Land Cover	\\
    \hline
CH-Dav	&	Switzerland	&	Davos	&	ENF	\\
DE-Tha	&	Germany	&	Tharandt	&	ENF	\\
FI-Hyy	&	Finland	&	Hyytiala	&	ENF	\\
IT-Lav	&	Italy	&	Lavarone	&	ENF	\\
RU-Fyo	&	Russia	&	Fyodorovskoye	&	ENF	\\
BE-Vie	&	Belgium	&	Vielsalm	&	MF	\\
CH-Lae	&	Switzerland	&	Laegern	&	MF	\\
DE-Hai	&	Germany	&	Hainich	&	DBF	\\
DK-Sor	&	Denmark	&	Soroe	&	DBF	
    \end{tabular}
    \end{ruledtabular}
    \label{table:sites}
    \end{table}

    \subsection{Methods: Time series analysis with dimension reduction techniques}
    \label{sec:methods}
    
    Dimension reduction techniques applied to time series data offer data-driven pattern extraction without \emph{a priori} assumptions about the nature of the time series itself.\cite{ghil2002advanced,2003KantzThomasSchreiber_Book,budivsic2012applied} Such spectral analysis approaches can yield interpretable features\cite{ghil2002advanced,giannakis2012nonlinear}. This study focuses on two dimension reduction methods: Singular Spectrum Analysis (SSA)\cite{broomhead1986extracting,ghil2002advanced} and Nonlinear Laplacian Spectral Analysis (NLSA)\cite{giannakis2012nonlinear}, which is based on Diffusion Maps.\cite{coifman2006diffusion} 
    Both methods, SSA and NLSA, optimize different projection metrics. SSA decomposes the covariance kernel, and minimizes the squared reconstruction error of all individual dimensions.\cite{van1999statistical} In contrast, NLSA aims to preserve the geometrical features associated with the discretized Riemannian metric\cite{berry2013timescale} by decomposing the diffusion kernel.    
    % different principles by which these features are extracted
    %To showcase the potential, strengths and challenges of trying to differentiate the seasonal cycle from the residual information with high precision, we compare the linear SSA and nonlinear NLSA across seven variables at nine measurement locations.
    
    We propose a multi-step approach using dimension reduction (SSA and NLSA) to detect harmonic oscillations in time series of land-atmosphere flux measurements, from which the seasonal cycle is constructed, Figure \ref{fig:workflow}:
    \begin{enumerate}
        \item The time series $x$ of length  $N=\qty{14}{\year}$, $x \in \R^N$ (referred to as signal), is embedded into a higher dimensional space employing Takens' delay coordinates. This yields a data matrix $X \in \R^{W,P}$, where $W = \qty{7}{\year}$ is the embedding window and $P$ is given by $P=N-W-1$.\cite{sauer1991embedology,takens1981detecting} The data matrix $X$ is standardized to zero mean and unit standard deviation per row.\cite{pearson1901lines} Centralizing the time series allows us to assign the individual components a rank order directly linked to the variance covered in the original time series (dimension reduction spectrum). 
        \item The dimension reduction algorithms computes dominant orthogonal vectors referred to as modes, which yield a low-dimensional representation of the embedded data $X$. This technique uses temporal dependencies to extract characteristic temporal patterns in the time series. 
        \item The modes representing pure harmonic oscillations are detected based on their spectral content using Fast Fourier Transform\cite{cooley1965algorithm} (FFT). As these modes have a physical meaning for the seasonal cycle, no grouping procedure typical for SSA is used.
        %. This aims to identify the individual modes that carry information solely related to the seasonal cycle. These modes represent the underlying periodic components of the time series and are filtered 
        Each harmonic cycle is comprised of mode pairs, i.e. the harmonic oscillation and its phase-shifted counterpart, as cycles have two degrees of freedom.
        %Oscillation with mixed frequencies are excluded from the reconstruction. 
        \item Finally, the seasonal cycle is constructed only from the harmonics and transformed into the original time domain, using a projection and reversing the centralization. 
        %(Appendix~\ref{app:centralization}).
    \end{enumerate}
    The mathematical details of these methods, along with the rationale for parameter choices, are provided in Appendix~\ref{app:methodsAll}.
    
    The numerical computations and figure generation were implemented in  Julia\cite{bezanson2017julia} and used the libraries \verb|MultivariateStats.jl|\cite{multivariatestats} and \verb|ManifoldLearning.jl|.\cite{manifoldlearning} The code to reproduce our results and figures is freely available.\cite{dimensionality} 
    
    To demonstrate the challenges posed by irregularity due to compromised data and high-frequency variability, we employed four different data characterization metrics, mainly in \S~\ref{sec:harmonicsVarSites}. We introduce these measures in an intuitive way here and give details in Appendix~\ref{app:dataChar}. We quantify \textbf{regularity} by normalizing the sum of the relative power of the first five harmonic oscillations extracted from the time series with FFT. We analyze both the unfiltered and the low-pass filtered signal, i.e. with frequencies $f \leq \qty{6}{\year^{-1}}$ and $\qty{4}{\year^{-1}}$, respectively. In addition, we quantify the \textbf{high-frequency variability} by computing the sum of the relative power sum of the signals high frequency content with $f > \qty{6}{\year^{-1}}$.
    To approximate the complexity, i.e. the randomness content, in our time series from an information theoretic point of view, we compute the \textbf{sample entropy}\cite{richman2000physiological} with standard choices such as delay $m=2$ and distance $r = \qty{0.2}{std}$.
    %, see Appendix~\ref{app:dataChar} for more details.
    We also account for artifacts resulting from data quality using \textbf{quality flags (QF)} at discrete time points. They indicate the confidence in the  measurement quality as well as in the data gap filling process in \%, i.e. a gray scale, where no confidence is indicated by black.
    We highlight irregular windows, which span QFs over at least $\unit{\year}/2$, in \S~\ref{sec:harmonicsVarSites}.

%#########################
%   S E C T I O N : Results
%########################

\section{Results}
\label{sec:results}

    \begin{figure*}[htb!]
        \includegraphics[width=\linewidth]{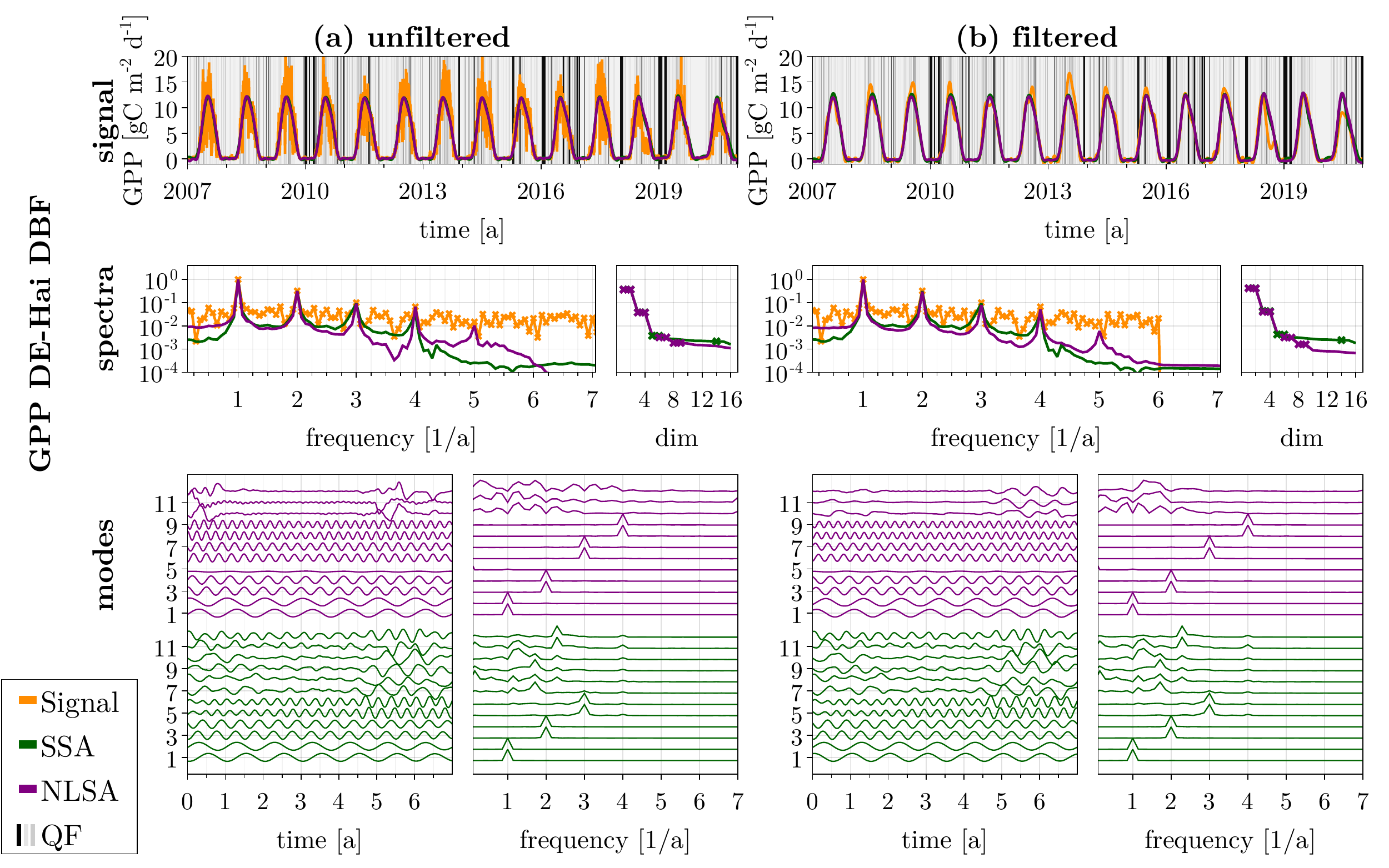}
        \caption{
        Analysis of a regular time series with isolated quality flags: gross primary production (GPP) fluxes of a deciduous broadleaf forest (DBF) in central Germany (Hainich forest site, DE-Hai DBF). \\
        Panels: Analysis of an unfiltered (a) and filtered (b) measurement signal. Each panel column illustrates the following analysis parts. 
        Panel \emph{signal}: signal and the corresponding constructed seasonal cycle with SSA or NLSA, if harmonics are detected. The quality flag (QF) is indicated by a gray scale, i.e. lowest signal quality in black.
        Panel \emph{spectra}: relative FFT power spectrum (left) of the signal and the SSA and NLSA constructed seasonal cycle, and the dimension reduction spectrum (right).
        Panel \emph{modes}: shapes of the 12 most dominant modes (left) and their corresponding spectra (right). \\
        Main result: The seasonal cycle is represented by a combination of harmonics of up to fifth order, i.e. SSA detects four harmonics and NLSA five harmonics.
        % old text 
        % (a) 4 (SSA) and 5 (NLSA) harmonic orders identified and accurate seasonal cycle for a relatively regular GPP signal. (b) Filtering (i.e. omitting all frequencies $>\SI{6}{a^{-1}}$) not affecting the seasonal cycle, only minor changes in decomposition spectrum and residual modes.
        }
        \label{fig:exRegular}
    \end{figure*}

%overview, dont kick down the door
The seasonal cycle in the extratropics is the oscillatory component of solar irradiation and closely linked climatic variables such as temperature with a fundamental frequency of $1 \unit{\year^{-1}}$. Land-atmosphere fluxes partially inherit these oscillatory components, but likewise represent more complex dynamics such as moisture availability, which is not necessarily strictly seasonal.  Hence, these ecologically mediated variables are inaccurately approximated by a first order sinusoidal function. Thus, in order to extract an ecologically more meaningful seasonal cycle, higher order harmonics are required.

We highlight the challenges of this approach related to measurement characteristics and irregularities by analyzing of four time series:
\begin{itemize}
    \item Case~1: regular time series 
    (Fig.~\ref{fig:exRegular})
    \item Case~2: time series with high frequency variability 
    (Fig.~\ref{fig:exResolved})
    \item Case~3: time series with broadband variability (Fig.~\ref{fig:exUnresolved})
    \item Case~4: time series with amplitude change
    (Fig.~\ref{fig:exDeficient})
\end{itemize}
Last, we compare the method's extraction performance across variables and measurement sites with different land cover types.

%we present the ecosystem  overview, which summarizes the analysis across all sites and variables.
%We considered the quality flag associated with the measurement time series and a FFT-based lowpass filter to omit high frequency variability.

%explanation what figures are we looking at -> DO NOT NEED THIS [KM]
%In the following we present four examples of vegetation flux time series and their lowpass-filtered versions. A worse confidence in measurement or gap filling in the quality flags (\%) is indicated by more pronounced grey scale vertical marks. The time series of the measurements are displayed together with the resulting computed seasonal cycle. Additionally, the relative harmonic power in the range of $\{0,\hdots,6\} \unit{\year^{-1}}$ by FFT highlights their respective spectral properties. The dimension reduction spectrum depicts the distribution of the dominant first projection values and offers insight into the differentiation and overlap of the seasonal cycle and residual information.
%Finally, the main modes of the two dimension reduction methods are displayed together with their spectral content in the range of $\{0,\hdots,7\} \unit{\year^{-1}}$.
%They were analyzed individually to reconstruct the seasonal cycle, consisting of the fundamental frequency and subsequent higher harmonic orders. An overview across all analyzed ecosystem is shown in Fig.~\ref{fig:harmonics}.

    \subsection{Case 1: Regular time series}
    \label{sec:case1}
    % \begin{figure*}[htb!]
    %     \includegraphics[width=\linewidth]{./pictures/ex_regular.pdf}
    %     \caption{
    %     Analysis of a regular time series with isolated quality flags: gross primary production (GPP) fluxes of a deciduous broadleaf forest (DBF) in central Germany (Hainich forest site, DE-Hai DBF). \\
    %     Panels: We analyze the unfiltered (a) and the filtered (b) measurement signal. Each panel column illustrates the following analysis parts. 
    %     Panel \emph{signal}: signal and the corresponding constructed seasonal cycle with SSA or NLSA, if harmonics are detected. The quality flag (QF) is indicated by a gray scale, i.e. lowest signal quality in black.
    %     Panel \emph{spectra}: relative FFT power spectrum (left) of the signal and the SSA and NLSA constructed seasonal cycle, and the dimension reduction spectrum (right).
    %     Panel \emph{modes}: shapes of the 12 most dominant modes (left) and their corresponding spectra (right). \\
    %     Main result: The seasonal cycle is represented by the harmonics of up to fifth order, i.e. SSA detects four harmonics and NLSA five harmonics.
    %     % old text 
    %     % (a) 4 (SSA) and 5 (NLSA) harmonic orders identified and accurate seasonal cycle for a relatively regular GPP signal. (b) Filtering (i.e. omitting all frequencies $>\SI{6}{a^{-1}}$) not affecting the seasonal cycle, only minor changes in decomposition spectrum and residual modes.
    %     }
    %     \label{fig:exRegular}
    % \end{figure*}

    We analyze GPP fluxes, i.e. the amount of assimilated CO$_2$, in a deciduous broadleaf forest  ecosystem in central Germany, Fig.~\ref{fig:exRegular}. These measurements display quality flags, Fig.~\ref{fig:exRegular}(a, signal), indicating a lower confidence in the respective measurement or gap filling method. The lowest quality flags (black) occur mostly in the winter months, where GPP is low and varies little. Despite such flags, our approach is capable of detecting four harmonics with SSA and five with NLSA, which indicates the signals regularity Fig.~\ref{fig:exRegular}(a, spectrum). The detection of harmonic oscillations is better with NLSA, as SSA exhibits mixed frequency signals for the third harmonic, e.g. for the mode corresponding to $f=\qty{3}{\per \year}$ we also observe a small peak at $f= \qty{4}{\per \year}$. Although the FFT power spectrum indicates peaks at the harmonic frequencies, they are less pronounced and could be overlooked. 
    
    The seasonal cycle is constructed from the harmonic modes and transformed into the time domain, Fig.~\ref{fig:exRegular}(a, signal). The seasonal cycles constructed from SSA and NLSA modes, respectively, are very similar and mostly overlap. They do not represent extra-seasonal events, such as the heatwave and drought in 2018, which leads to decreased GPP, i.e. decreased vegetation productivity, in autumn. 
    
    Although filtering the time series reduces the signal's variability greatly, Figure \ref{fig:exRegular}(b, signal), the SSA and NLSA results change little, especially the number of detected harmonics. The filtered signal and the constructed seasonal cycle overlap mostly. Note that, the filtered signal does retain the impact of extra-seasonal events, such as the heatwave and drought in 2018. A similar case is shown in Appendix~\ref{app:suppFig} Fig.~\ref{fig:exAppendix1}.

    % OLD TEXT from Leonard
    % Figure \ref{fig:exRegular}(a) displays Gross Primary Production (GPP) fluxes - the amount of assimilated CO$_2$ - in a deciduous broadleaf forest (DBF) ecosystem in Germany (the Hainich national park measurement site DE-Hai). In many instances of short duration quality flags indicate lower confidence in the respective measurement or gap filling. Appearing mostly in the winter months of very low GPP variability, the signal's regularity is not hampered, i.e. the ability to extract higher order harmonics. The time series appears very regular and the higher harmonic orders were obtained that describe the seasonal cycle both by SSA (4 harmonics) and NLSA (5 harmonics). The dimension reduction spectrum of the dimension reduction has a clear slope-break around the third harmonic. Each frequency is resolved by two corresponding linearly independent modes. The modes from third order and up of SSA display an information overlap with the fourth order and lower-frequency residual information.

    % Figure \ref{fig:exRegular}(b) displays the lowpassed signal of Figure \ref{fig:exRegular}(a). In (a), the filtered signal showed slight amplitude modulations of the summer maxima throughout the years, that were not visible among the unfiltered high frequency variability. The resulting seasonal cycle of SSA does not display change after omitting the high frequency variability. For NLSA, the modes directly after the harmonic structure now showed marginally more coherent low-frequency spectral content and a minor change in the dimension reduction spectrum slope.
    
    \subsection{Case 2: Time series with high frequency variability}
    \label{sec:case2}
    
    \begin{figure*}[ht!]
        \includegraphics[width=\linewidth]{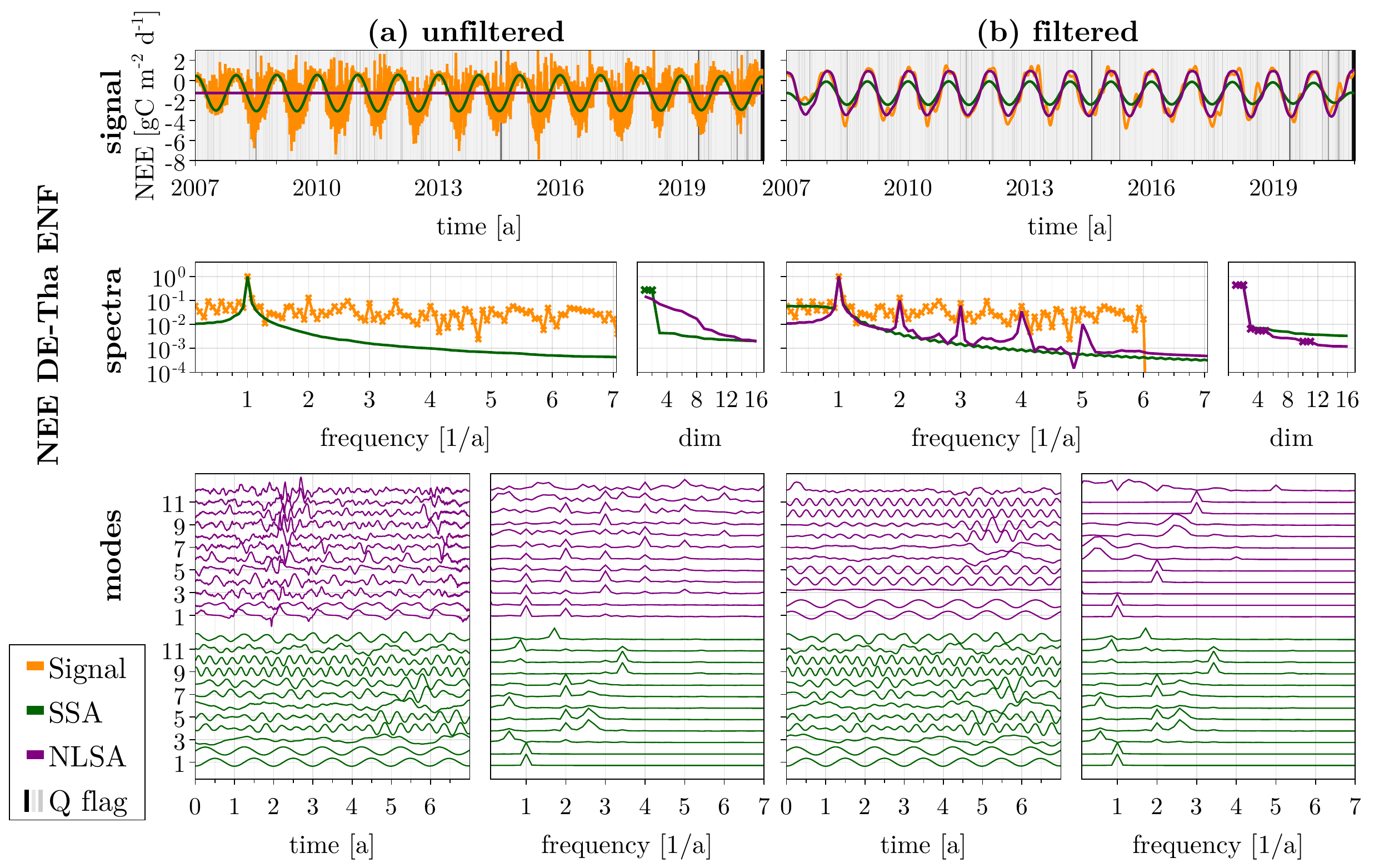}
        \caption{
        Analysis of a time series with high frequency variability and isolated quality flags: net ecosystem exchange (NEE) fluxes of a evergreen needleleaf forest (ENF) in Eastern Germany (Tharandt forest site, DE-Tha ENF). \\
        Panels: We analyze the unfiltered (a) and the filtered (b) measurement signal. Each panel column illustrates the following analysis parts. 
        Panel \emph{signal}: signal and the corresponding constructed seasonal cycle with SSA or NLSA, if harmonics are detected. The quality flag (QF) is indicated by a gray scale, i.e. lowest signal quality in black.
        Panel \emph{spectra}: relative FFT power spectrum (left) of the signal and the SSA and NLSA constructed seasonal cycle, and the dimension reduction spectrum (right).
        Panel \emph{modes}: shapes of the 12 most dominant modes (left) and their corresponding spectra (right). \\
        Main result: The harmonics of up to fifth order are only detected after filtering, i.e. SSA detects only the fundamental oscillation while NLSA detects five harmonics. The amplitudes of the seasonal cycles constructed with SSA and NLSA, respectively, differ. 
        % OLD caption        
        %NEE measurements from DE-Tha site (ENF). Irregular time series with strong high frequency variability. Lowpass-filtering enables detection of seasonal cycle by NLSA. Identical layout as Fig.~\ref{fig:exRegular}. (a) Seasonal cycle resolved only with fundamental frequency by SSA, not at all by NLSA. Individual NLSA modes show overlap of harmonic frequencies. (b) After filtering, NLSA detects 5 orders of harmonics. No visible change for SSA. Second order subharmonics are visible across different modes, i.e. $f = \{1/2,3/2,5/2\}\unit{\per\year}$.
        }
        \label{fig:exResolved}
    \end{figure*}

    We analyze NEE, the net exchange of carbon including both CO$_2$ uptake by the ecosystem via photosynthesis and release of CO$_2$ via respiration, at an evergreen needleleaf forest site in Eastern Germany, Fig.~\ref{fig:exResolved}. The signal is affected only by minor quality flags (light gray). However, it shows a high variability across time, which affects the detection of harmonics modes.

    With SSA we extract only the fundamental mode, irrespective of signal filtering, Fig.~\ref{fig:exResolved}(a \& b, spectra). The other modes exhibit oscillations of mixed harmonic and non-harmonic frequencies, indicating that SSA cannot separate these along the characteristic directions, i.e. the maximum data variance, Fig.~\ref{fig:exResolved}(a \& b, modes, green).

    NLSA could not separate the harmonic oscillation into single modes before filtering out the high frequency signal content, Fig.~\ref{fig:exResolved}(a, modes, purple). The FFT spectrum reveals, that the power of the high frequency oscillations is comparably high to that of the other oscillations, Fig.~\ref{fig:exResolved}(a, spectra, left). Unlike SSA, the dimension reduction spectrum does not show a slope break, Fig.~\ref{fig:exResolved}(a, spectra, right). Only after filtering did NLSA detect five purely harmonic modes, Fig.~\ref{fig:exResolved}(b, spectra, right). The detected modes are attributed a lower power than in the previous cases, i.e. to dimensions greater than 12, Fig.~\ref{fig:exResolved}(b, modes, purple)

    Consequently, the seasonal cycles constructed from the harmonic oscillations extracted with SSA and NLSA, respectively, exhibit a clear difference in amplitude, Fig.~\ref{fig:exResolved}(b, signal).
    A similar case is shown in Appendix~\ref{app:suppFig} Fig.~\ref{fig:exAppendix2}.

    \subsection{Case 3: Time series with broadband variability}
    \label{sec:case3}
    
    \begin{figure*}[ht!]
        \includegraphics[width=\linewidth]{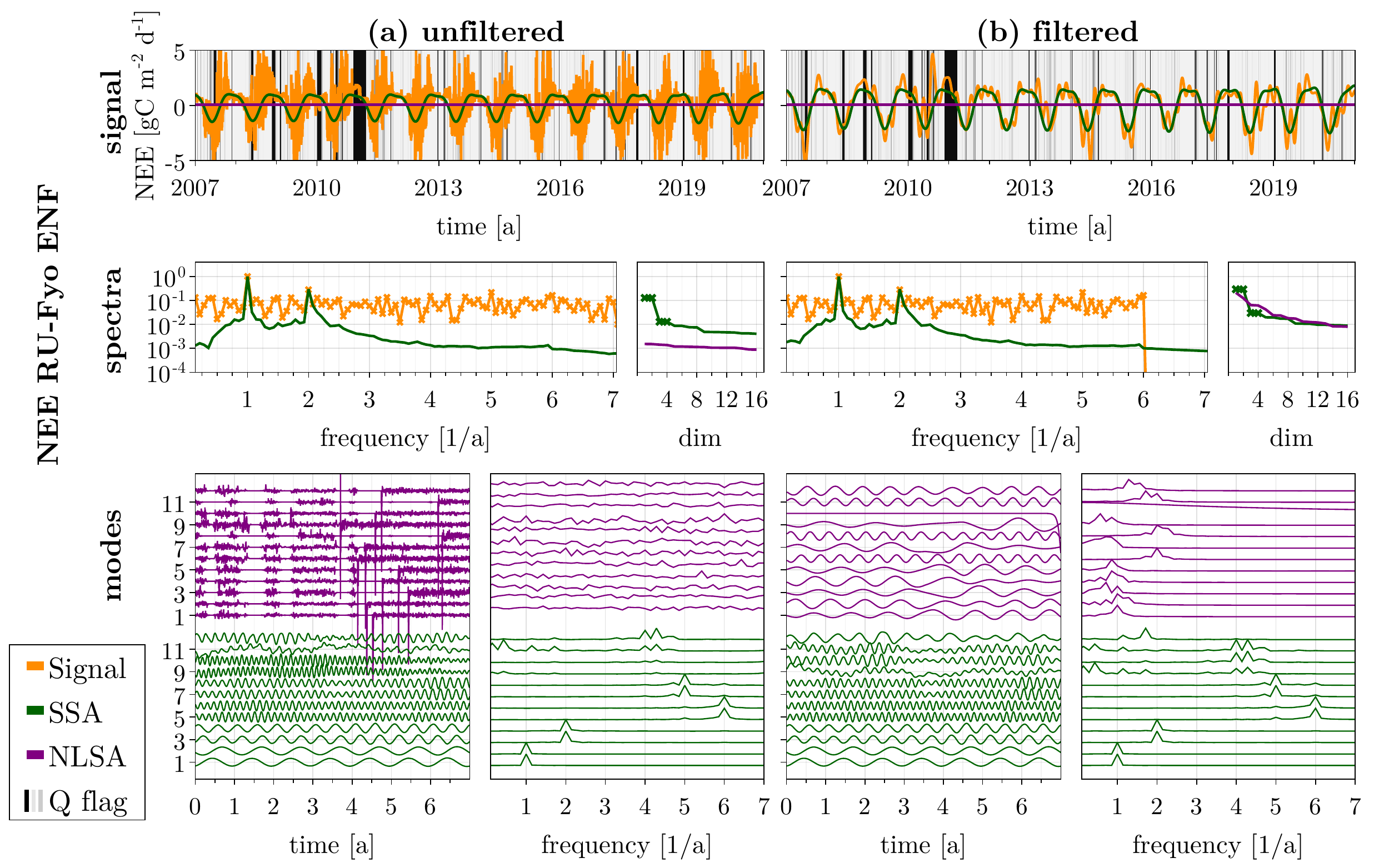}
        \caption{
        Analysis of a time series with broadband variability and poor quality flags over extended periods of time: net ecosystem exchange (NEE) fluxes of an evergreen needleleaf forest (ENF) in western Russia (Central Forest Biosphere Reserve at Fyodorovskoe site, RU-Fyo ENF). 
        \\
        Panels: Analysis of an unfiltered (a) and filtered (b) measurement signal. Each panel column illustrates the following analysis parts. 
        Panel \emph{signal}: signal and the corresponding constructed seasonal cycle with SSA or NLSA, if harmonics are detected. The quality flag (QF) is indicated by a gray scale, i.e. lowest signal quality in black.
        Panel \emph{spectra}: relative FFT power spectrum (left) of the signal and the SSA and NLSA constructed seasonal cycle, and the dimension reduction spectrum (right).
        Panel \emph{modes}: shapes of the 12 most dominant modes (left) and their corresponding spectra (right). \\
        Main result: After filtering, the time series exhibits a level of variability, which still prohibits the detection of any harmonics with NLSA. In contrast, SSA always detects the fundamental and the second harmonic oscillation. 
        % OLD caption
        %NEE measurements from RU-Fyo site (ENF). Irregular time series with high frequency variability and severe data quality issue. Identical layout as Fig. \ref{fig:exRegular}. The quality flags at 2011 mark a jump in the measurement data. (a) SSA detects fundamental and second order harmonic. NLSA modes erratic shaped and without distinct spectral content. NLSA decomposition spectrum almost flat, showing no distinction between information. (b) After filtering, NLSA decomposition spectrum increased slope and modes become more shaped like harmonic structure - not detected still. Minor changes in SSA residual mode appearance.
        }
        \label{fig:exUnresolved}
    \end{figure*}

    We analyze NEE, the net exchange of carbon including both CO$_2$ uptake by the ecosystem via photosynthesis and release of CO$_2$ via respiration, at an evergreen needleleaf forest site in Western Russia, Fig.~\ref{fig:exUnresolved}. This time series is affected by poor quality flags over prolonged time periods, Fig.~\ref{fig:exUnresolved}(a \& b, signal). It is characterized by a strong broadband variability, obviously unaffected by the low-pass filter, Fig.~\ref{fig:exUnresolved}(a \& b, spectra, left).

    SSA successfully extracts the fundamental and the second harmonic oscillations from both the unfiltered and filtered signal, Fig.~\ref{fig:exUnresolved}(a \& b, spectra, left). However, while the fourth, fifth, and sixth harmonics dominate modes 5-8 and 11,12, SSA is unable to separate these oscillations from the nonharmonic ones, Fig.~\ref{fig:exUnresolved}(a \& b, modes). Note, that the third harmonic appears to exhibit a low variance, and is not detected within the first 16 dimensions. Moreover, the constructed seasonal cycle does not capture the signal‘s broadband frequency variability, Fig.~\ref{fig:exUnresolved}(b, signal).

    NLSA fails to extract any pure oscillations. This is evidenced by the dimension reduction spectrum, which does not show a slope break (Fig.~\ref{fig:exUnresolved}(a, spectra, right)), and by the noisy and shapeless modes, Fig.~\ref{fig:exUnresolved}(a, modes). After filtering, the modes are primarily characterized by the fundamental and second harmonic oscillations (Fig.~\ref{fig:exUnresolved}(b, modes)), yet they still represent mixed oscillations and cannot  be separated.
    A similar case is shown in Appendix~\ref{app:suppFig} Fig.~\ref{fig:exAppendix3}.
    
    %Figure \ref{fig:exUnresolved} displays NEE fluxes in an evergreen needleleaf forest at the RU-Fyo measurement site. The time series is characterized by a constant large variability and a major quality flag indication in 2011. There, the measurements are constant before dropping abruptly. SSA computed two harmonic orders, yielding a characteristic NEE seasonal cycle. The spectra of the modes reveal a gap in the harmonic structure at third and fourth order. The computed NLSA dimension reduction spectrum has an almost linearly decreasing slope, not identifying any part of the seasonal cycle.   The individual NLSA modes appear with noise-like variability and showed no distinct spectral content.
    %
    %In Figure \ref{fig:exUnresolved}(b), the filtered time series appear irregular still. The computed seasonal cycle by SSA did not change in mode composition. The NLSA modes on the other hand reveal an oscillatory shape, while their spectra did not qualify them as purely harmonic. The spectral content clusters around first and second frequency.
    
    \subsection{Case 4: Time series with amplitude change}
    \label{sec:case4}
    
    \begin{figure*}[ht!]
        \includegraphics[width=\linewidth]{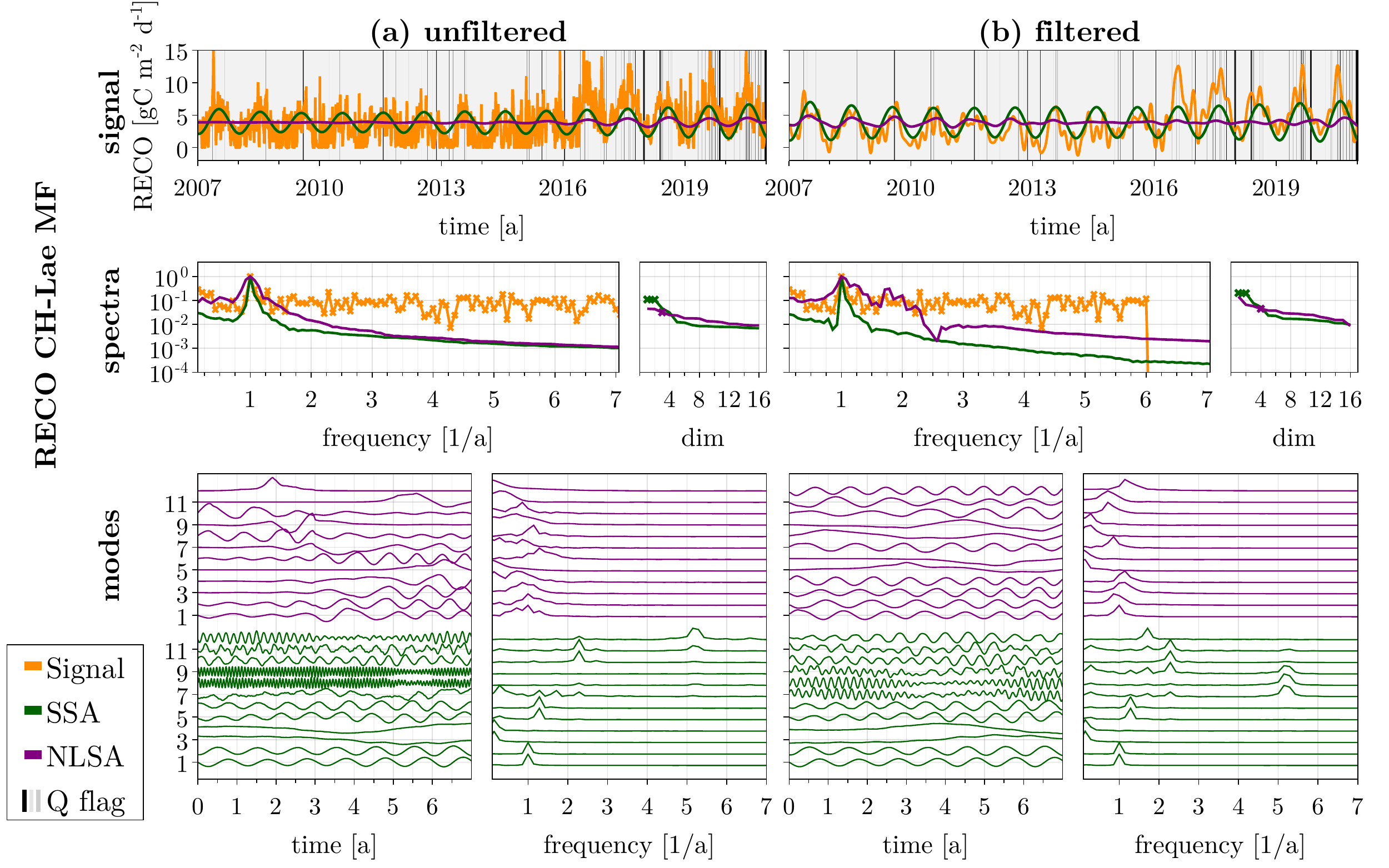}
        \caption{
        Analysis of a time series with amplitude change in addition to broadband variability and isolated quality flags: ecosystem respiration ($\reco$) of a mixed forest (MF) in Northern Switzerland (Laegeren mountain site, CH-Lae MF). 
        \\
        Panels: Analysis of an unfiltered (a) and filtered (b) measurement signal. Each panel column illustrates the following analysis parts.
        Panel \emph{signal}: signal and the corresponding constructed seasonal cycle with SSA or NLSA, if harmonics are detected. The quality flag (QF) is indicated by a gray scale, i.e. lowest signal quality in black.
        Panel \emph{spectra}: relative FFT power spectrum (left) of the signal and the SSA and NLSA constructed seasonal cycle, and the dimension reduction spectrum (right).
        Panel \emph{modes}: shapes of the 12 most dominant modes (left) and their corresponding spectra (right). \\
        Main result: The low-pass filter shows the amplitude change after 2016. SSA extracts the fundamental oscillation irrespective of noise. NLSA extracts a single mode (not a pair!) due to the error bounds within our algorithm. This clearly does not represent a pure oscillation.
        % OLD caption
        %$\reco$ measurements from CH-Lae site (MF). Irregular signal with high frequency variability and trend/drift/transition: much stronger amplitudes, more pronounced seasonal cycle after 2016. Example for deficient harmonic structure. 
        }
        \label{fig:exDeficient}
    \end{figure*}

    We analyze $\reco$, the total ecosystem respiration, which is the sum of the carbon loss during plant respiration and decomposition, at a mixed forest site in Northern Switzerland, Fig.~\ref{fig:exDeficient}. This time series, with few quality flags, exhibits nonstationary behavior, i.e. undergoes a fast amplitude change in 2016. The low-pass filter accentuates this transition even clearer, Fig.~\ref{fig:exDeficient}(a \& b, signal). The signal exhibits a strong broadband variability, Fig.~\ref{fig:exDeficient}(a) \& (b, spectra, left) as in Case~3.

    Despite the nonstationary nature of the data, SSA extracts the fundamental oscillation. The flat slope break in the dimension reduction spectrum, Fig.~\ref{fig:exDeficient}(a \& b, spectral, right), suggests that many characteristic directions are required to represent the signal. Consequently, the seasonal cycle, constructed only from the fundamental oscillation, fails to capture the amplitude change or the variability, Fig.~\ref{fig:exDeficient}(a \& b, signal).

    As shown in Case~3, NLSA struggles to extract harmonic oscillations, when the signal is affected by strong broadband variability. Here, NLSA detects a single mode with the fundamental frequency, however, this is primarily a result of the algorithm‘s built-in error bounds. The corresponding FFT spectrum shows that this is not a pure oscillation. The detected modes are primarily mixed oscillations with frequencies below 1.5. 
    A similar case is shown in Appendix~\ref{app:suppFig} Fig.~\ref{fig:exAppendix3}.
    
    %Other site variables show analogous behavior, see Figs.~\ref{fig:exAppendix1}--\ref{fig:exAppendix3} in Appendix~\ref{app:suppFig}.
    %In Figure \ref{fig:exDeficient}, an example of a challenge for the extraction of the seasonal cycle by SSA/NLSA is displayed.
    %
    %- time series visually exhibits nonstationary behavior (different amplitude after 2016) (is this still an anomaly?) 
    % - ssa makes out clear fundamental, all other modes have overlap in spectral content
    % - nlsa modes all cluster around fundamental (modulated modes with broad spectral peaks)
    % - dimensionality reduction spectrum is missing the typical slopebreak behavior, possibly indicating issues in the differentiation between seasonal and residual information
    % - filtering does not resolve the missing resolution as it did in Fig. \ref{fig:exUnresolved}, i.e. the issue is not (only) the signal/noise ratio
    % - seasonal cycle not identified by either method SSA or NLSA (due to nonstationary behavior)
    % - interpretation (results): closest explanation is missing samples (9 years small amplitude, 5 years larger amplitude), but a shift in amplitude might be a no-go for dimensionality reduction in this fashion after all
    
    % Similar results are obtained by climatic variables, see 
    % Figures \ref{fig:AppendixunfilteredSeason} and 
    % \ref{fig:AppendixfilteredSeason} in Appendix \ref{app:suppFig}.

    \subsection{Harmonic extraction across variables and measurement sites}
    \label{sec:harmonicsVarSites}
    \label{sec:allSites}
    
    \begin{figure}[htb!]
        \centering
        \includegraphics[width=\linewidth]{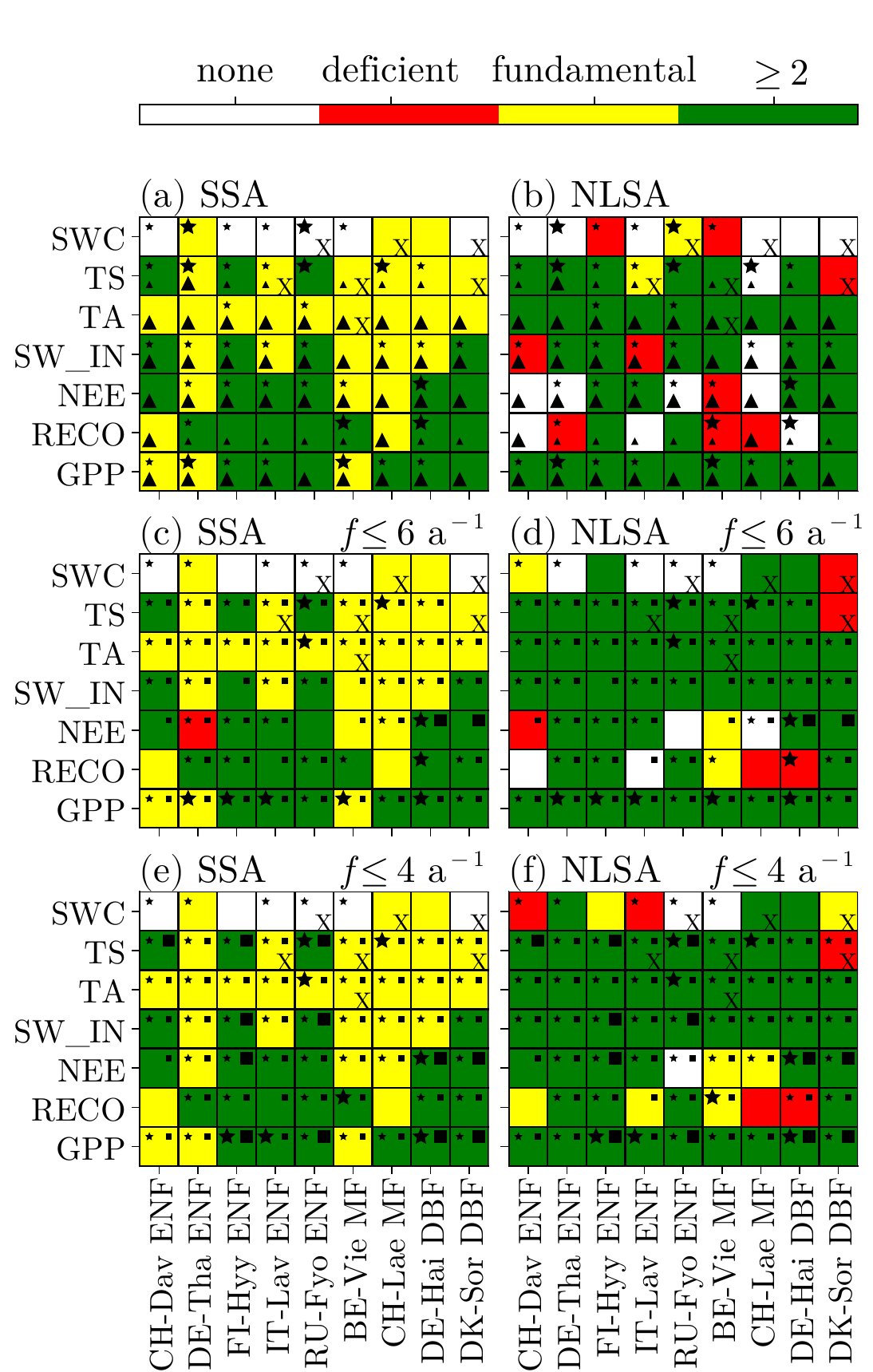}
        \caption{
        Harmonic oscillation extraction per time series across variables and measurement sites with SSA (a,c,e) and NLSA (b,d,f).
        Sites (Table~\ref{table:sites}) with land cover types ($x$-axis): evergreen needleleaf forest (ENF), mixed forest (MF), deciduous broadleaf forest (DBF).
        Variables ($y$-axis): soil water content (SWC), soil temperature (TS), air temperature (TA), Short wave radiation (SW\_IN), net ecosystem exchange (NEE), ecosystem respiration ($\reco$), gross primary productivity (GPP).
        Colors indicate whether harmonic oscillations have been detected in the modes: 
        white/none = no pure oscillations detected; 
        red/deficient = fundamental not resolved; 
        yellow/fundamental = only fundamental oscillation extracted (a mode pair);
        green/$\geq2$: at least the fundamental and the second harmonic extracted. 
        Filters: 
        unfiltered (a,b), 
        low-pass filter with $f = \qty{6}{\per \year}$ (c,d) 
        and with $f = \qty{4}{\per \year}$ (e,f). 
        Markers for time series properties:
        X = persistent quality flags over at least $\unit{a}/2$,
        $\blacktriangle$ = high frequency variability (large, low, none),
        $\blacksquare$ = regularity (high, low, none), and 
        {\Large $\star$} = sample entropy (high, low, no randomness).
        % low  sample entropy = higher regularity 
        Main result: SSA fails to detect any pure oscillations (none \& deficient) in at most seven time series (c), compared to NLSA which fails in at most 24 time series (b). But, NLSA extracts higher order harmonics ($\geq2$) in more time series (at most 48 in (f)) than SSA (at most 26 in (a,c,e)). The low pass filter improves extraction with NLSA, but not with SSA. 
        % 1) quality flags of at least $\text{a}/2$ length (X),
        % 2) strength of high frequency variability, marked by large, small or no $\blacktriangle$,
        % 3) sample entropy, marked by large, small or no $\star$,
        % 4) regularity, marked by large, small or no $\blacksquare$.
        % OLD TEXT
        % classification of seasonal cycle quality across multiple variables and ecosystems.
        % white:none (fundamental not detected);
        % red:deficient (fundamental not resolved accurately);
        % yellow:fundamental (only sinusoidal 1a);
        % green:at least second harmonic order detected.
        % Left side: SSA, Right side: NLSA,
        % top row (a,b) : unfiltered,
        % mid row (c,d) : filter as in examples at $f_l = \SI{6}{\year^{-1}}$,
        % bottom row (e,f) strong filter at $\SI{4}{\year^{-1}}$.
        % Time series properties are marked as follows:
        % 1) quality flag artifacts of at least $\text{a}/2$ length (X for existence),
        % 2) strength of high frequency variability, marked by large, small or no $\blacktriangle$,
        % 3) sample entropy, marked by large, small or no $\star$,
        % 4) regularity, marked by large, small or no $\blacksquare$.
        % No effect on SSA by filtering, green correlates highly with regularity,
        % many cases of only fundamental.
        % NLSA is helped by $\SI{6}{\year^{-1}}$ filter,
        % no great improvement by further lowpassing with $\SI{4}{\year^{-1}}$.
        }
        \label{fig:harmonics}
    \end{figure}
    
    % What is the aim of this section and what to we explore here
    In this section, we examine the capabilities of SSA and NLSA to extract harmonic oscillations across seven variables and nine measurement sites with various land cover types (ENF, MF, DBF), Fig.~\ref{fig:harmonics}. Carbon fluxes, specifically GPP, NEE, and $\reco$
    %Gross Primary Production (GPP), Net Ecosystem Exchange (NEE), and ecosystem respiration ($\reco$),
    , exhibit site-specific dynamics and variability due to local ecosystem conditions such as soil water content and soil temperature, in addition to local ecosystem drivers such as air temperature and short wave radiation.

    % What do we show here
    In Fig.~\ref{fig:harmonics} We analyze both the unfiltered and filtered signals, employing bandpass filters with frequencies $f \leq \qty{6}{\per \year}$ and $f \leq \qty{4}{\per \year}$. To simplify the visualization, we categorize the detection of harmonic oscillations into four distinct categories:
    \begin{itemize}
        \item no pure oscillations are detected, as seen in Case~3 with NLSA (none).
        \item only one mode of a pair is detected, typically the mode with the fundamental frequency, as observed in Case~4 with NLSA (deficient).
        \item only the fundamental frequency is detected, exemplified by Case~2 with SSA (fundamental).
        \item multiple higher-order harmonics are detected by mode pairs, such as in Case~1 with both SSA and NLSA ($\geq 2$).
    \end{itemize}
    The specific number of harmonics is shown in Appendix~\ref{sec:harmonicsNum}.

    % overall result of the figure: number of harmonics extracted
    Overall, SSA fails to detect any oscillatory components in fewer time series than NLSA, see Fig.~\ref{fig:harmonics} (none and deficient categories) and Table~\ref{table:harmonics} (No H). However, when oscillations are detected, NLSA identifies higher-order harmonics in more cases than SSA (Table~\ref{table:harmonics}, H $\geq 2$), while SSA detects only the fundamental oscillation in almost half the cases (Table~\ref{table:harmonics}, H $=1$). The detection of harmonics using NLSA improves substantially when the signal is filtered.
    We find that the number of extracted harmonics is variable specific rather than site specific. For example, both SSA and NLSA have difficulty extracting harmonics from soil water content, as this is an aperiodic signal with high variability across frequency scales (Fig.~\ref{fig:exAppendix3}). 

    % data characterization measures
    To quantify this signal's variability we employ the following four measures: 
    persistent quality flags, i.e. QF over a continuous period of at least $\unit{a}/2$ (X),
    strength of high frequency variability ($\blacktriangle$), 
    regularity ($\blacksquare$), 
    and sample entropy ($\star$).

     % measure: quality flags
    Persistent quality flags (X) are found in seven instances. In more than half the instances these do not have an affect on the detection of harmonics. SSA fails to detect any oscillations in two instances, while NLSA fails in two to three depending on the filter. Filtering does not improve the detection substantially, indicating the irregularity is persistent and does not have a high frequency character. SSA is at best able to detect the fundamental, while NLSA detects higher order harmonics. 
    
    % measure: High-frequency variability
    High-frequency variability ($\blacktriangle$) is consistently observed across sites for variables such as air temperature, short wave radiation, and the carbon fluxes NEE and GPP; but not for soil water content, soil temperature, and ecosystem respiration, Fig.~\ref{fig:harmonics} (a,b). The low-pass filters effectively reduce this high-frequency content across all sites, Fig.~\ref{fig:harmonics} (c-f). As demonstrated in Cases~1--4, filtering enhances the detection of higher-order harmonics with NLSA but not SSA. The reduction in high-frequency content through filtering allows NLSA to more accurately identify and extract these harmonics, thereby improving its performance in analyzing oscillatory behavior.

    % measure: regularity
    The time series regularity detected with FFT emerges after filtering, Fig.~\ref{fig:harmonics} (a,b) vs (c-f). This indicates that high frequency content affects the detection of harmonics with FFT. SWC is the only variable with no regularity irrespective of filtering, indicating its aperiodic nature. High regularity in time series coincides with detection of higher order harmonics with SSA and NLSA, indicating that the three methods agree, especially in Fig.~\ref{fig:harmonics} (e,f). No regularity indicates that either SSA or NLSA or both fail to detect harmonics.

    % measure: sample entropy
    % low values: high self-similarity, regular
    % high values: random time series
    Sample entropy is a statistic that measures the regularity of a time series: low values represent regular time series, whereas high values indicate random behavior. Unlike the other fairly periodic variables, soil water content time series are usually aperiodic (Fig.~\ref{fig:exAppendix3}). Against our expectation, the sample entropy for SWC is at most low after filtering. Overall, time series with aperiodic behavior or amplitude change, did not have high entropy, Figs.~\ref{fig:harmonics}, \ref{fig:exDeficient}, and \ref{fig:exAppendix3}. We find that the sample entropy is affected by filtering, i.e. it can both decrease and increase, Fig.~\ref{fig:harmonics}.
    Moreover, higher order harmonics are detected in time series with high sample entropy in more than half the cases, Fig.~\ref{fig:harmonics}. This indicates a regular rather than random signal. 

    % nonstationary time series: aperiodic behaviour or amplitude change
    We aimed to use SWC time series for additional validation, by testing when SSA and NLSA fail to detect harmonics. We find that this failure is a better indicator for aperiodic or nonstationary time series (bandpass filtered with $f\leq \qty{6}{\per \year}$) than the sample entropy, Figs.~\ref{fig:harmonics}. We focus on the filtered time series here, as high-frequency content can affect NLSA's performance. 
    In addition, SSA and/or NLSA detects harmonic oscillations in SWC, despite its aperiodic nature, Figs.~\ref{fig:harmonics}. These seasonal cycles only loosely relate to the full signal, Fig.~\ref{fig:exAppendix3}.

    \begin{table}[h!]
    \caption{Percentage of sites, where SSA or NLSA detect (i) no pure oscillatory or deficient modes (no H), (ii) only the fundamental oscillation (H $=1$), or (iii) at least the first two harmonics (H $\geq2$). The time series are either unfiltered (filter none), or band-pass filtered with $f \leq \qty{6}{\per \year}$ or $f \leq \qty{4}{\per \year}$. This table summarizes the results from Fig.~\ref{fig:harmonics}.
    }
    \centering
    \begin{ruledtabular}
    \begin{tabular}{c l  r r r  r r r }
    & & \multicolumn{3}{ c  }{\textbf{SSA}} & \multicolumn{3}{ c }{\textbf{NLSA}} \\
    %\rowcolor{gray!30}[][]
    & \textit{filter} &	\textit{none}	&	\textit{$f\leq \qty{6}{\per \year}$}	&	\textit{$f\leq \qty{4}{\per \year}$	}&	\textit{none}	&	\textit{$f\leq \qty{6}{\per \year}$}	&	\textit{$f\leq \qty{4}{\per \year}$}	\\
    %harmonics & & & & & & \\
    \hline
    \parbox[t]{1mm}{\multirow{4}{*}{\rotatebox[origin=l]{90}{\% \textbf{Sites}}}}
&	no H &	9,5\%	&	11,1\%	&	9,5\%	&	38,1\%	&	20,6\%	&	12,7\%	\\
&	H $=1$	&	49,2\%	&	47,6\%	&	49,2\%	&	3,2\%	&	4,8\%	&	11,1\%	\\
&	H $ \geq2$	&	41,3\%	&	41,3\%	&	41,3\%	&	58,7\%	&	74,6\%	&	76,2\%	\\

    \end{tabular}
    \end{ruledtabular}
    \label{table:harmonics}
    \end{table}
         
%#########################
%   S E C T I O N : Discussion
%########################
\section{Discussion}
\label{sec:discuss}    

    \subsection{Seasonal cycle}
    \label{sec:seasonalCycle}
    % Seasonal cycle is comprised of harmonic oscillations
    In this study, we use data-driven approaches to uncover the inherent dynamics of the seasonal cycle of vegetation carbon processes. Our results show that the seasonal cycle is comprised of more than just the fundamental oscillation.
    We show that the seasonal cycle can be effectively approximated by a linear combination of harmonic modes of up to fifth order, Fig.~\ref{fig:exRegular}. These harmonics not only indicate the underlying complexity of vegetation process measurements and their drivers, but also emphasize the potential of data-driven approaches in capturing and unraveling intricate temporal patterns. 
    Although the computational principle of our two methods, SSA and NLSA, differ, they both reveal harmonic oscillations.
    Good agreement is found with the frequencies obtained with the FFT (of the entire signal), which uses sinusoidal basis functions, unlike SSA and NLSA. For example, other studies using FFT \cite{linscheid2021time} detect and construct the seasonal cycle from at most two harmonics.
    Our data-driven approach aligns with the historical approach of modeling the seasonal cycle in atmospheric carbon concentration measurements with sinusoidal functions\cite{Bacastow1985seasonal}, which are linear combinations of harmonics with frequencies $f=1, 2, \ldots, n \, \unit{\per \year}$. 
    %Our study also aligns with the historical approach of modeling the seasonal cycle in Earth system process measurements with sinusoidal functions, which are linear combinations of harmonics with frequencies $f=1, 2, \ldots, n$, such as the atmospheric carbon dioxide measurements \cite{Bacastow1985seasonal}. 

    % Number of harmonics, slope break
    Our findings demonstrate that the number of harmonic oscillations indicates the complexity of the seasonal cycle of carbon fluxes, e.g. Fig.~\ref{fig:exRegular}.
    Detecting an insufficient number of harmonics can lead to a poor approximation of the seasonal cycle. For example, when SSA detects only the fundamental component of the seasonal cycle the signals amplitude is misrepresented, Fig.~\ref{fig:exResolved}.
    Typically, when applying dimension reduction approaches to time series, modes are computed to maximize the amount of variance (related to the spectrum) they can explain, e.g. see spectrum-dimension plots in Fig.~\ref{fig:exRegular}. At some point, modes corresponding to higher dimensions contribute negligibly to explaining data variance. This point is referred to as the slope break in the spectrum-dimension plot. However, we show that for carbon fluxes, this is not a sufficient criterion to estimate a sufficient number of harmonic components. For instance, in Fig.~\ref{fig:exResolved}(b) the NLSA slope break occurs at dimension three, which only yields the fundamental oscillation, misrepresenting the seasonal cycle's amplitude (as shown by SSA seasonal cycle).

    % harmonics and data variance, detection order
    Our analysis also shows that the harmonic character of the signal does not necessarily align with its variance. In other words, the order of the modes dictated by the variance does not correspond to the order of the harmonic frequencies. 
    %For example, the SSA mode of frequency $f=4$ and the NLSA mode of frequency $f=5$ are not in the set of the first twelve modes, while the slope break is estimated to occur at dimension six, Fig.~\ref{fig:exRegular}(a). 
    For example, NLSA mode pairs $(1,2)$ correspond to the fundamental frequency $f=\qty{1}{\per \year}$, but higher order harmonic oscillations with frequencies $f=\qty{2}{\per \year}$ and $\qty{3}{\per \year}$ are found in mode pairs $(4,5)$ and $(10,11)$, respectively, Fig.~\ref{fig:exResolved}(b).
    Thus, a slope break at low dimensions may fail to detect higher order harmonics.
    
    % noise, variability 
    The detection of the seasonal cycle is affected by signal extra-seasonal anomalies such as extreme weather events, artifacts, and noise. The presence of noise can hinder the detection of harmonic oscillations, e.g. with NLSA as illustrated in Fig.~\ref{fig:exResolved}(a). However, effective noise filtering holds the potential to facilitate successful detection, as depicted in Fig.~\ref{fig:exResolved}(b). In the same case, we demonstrate that the intricacies inherent to the seasonal cycle can inadvertently evade detection (here by SSA), potentially leading to an underestimation of its amplitude despite filtering, Fig.~\ref{fig:exResolved}(b). Consequently, such an underestimation could affect any analysis relying on deseasonalization, which extracts the seasonal cycle to study variability.\cite{mahecha_characterizing_2007}

    \subsection{SSA vs. NLSA: comparing the linear and the nonlinear method}
    \label{sec:SSAvsNLSA}
    % differences maths principals of SSA and NLSA and relation to results: 
    The different computational principals of SSA and NLSA are reflected in their detection of harmonics. The graph-based algorithm NLSA detects dynamically important components, here higher order harmonics, Table~\ref{table:harmonics} and Fig.~\ref{fig:harmonics}. In contrast, SSA, which computes components by maximizing the time series' variance, mainly detects the dominant oscillations, i.e. often only the fundamental oscillation, Table~\ref{table:harmonics} and Fig.~\ref{fig:harmonics}. In such cases, the other modes contain mixed oscillation with harmonic characteristics, which SSA cannot separate. 
    Similar observations are made, for sea surface temperature model time series, where NLSA, but not SSA, detects dynamically important modes such as low frequency and intermittent modes.\cite{giannakis2012nonlinear} 
    Thus, our study highlights that land-atmosphere exchange of CO$_2$ involves complex processes where not all dynamical aspects are represented by components with a large variance.
    
    % more details on harmonics. -> I think it is too many details.
    % The linear method SSA has a lower failure rate than NLSA, Table~\ref{table:harmonics}.
    % ne of the main differences between the two methods is the detection success of harmonic oscillations. The linear method, SSA, successfully detects at least the fundamental oscillation more often than the nonlinear method, NLSA, see Fig.~\ref{fig:harmonics}. At the same time, when NLSA does detect harmonic oscillation, it captures more than SSA. An example is shown in Fig.~\ref{fig:exResolved}(b), where SSA detects the fundamental oscillation, while NLSA detects five harmonic oscillation pairs. Consequently, when an intricate harmonic structure in the signal is presumed, such as in vegetation carbon process measurements (GPP, NEE, $\reco$), then NLSA is better suited than SSA to uncover it. 
    % % When the signal is fairly regular with little noise, both methods uncover intricate harmonic structures, Fig.~\ref{fig:exRegular}.

    %Computational success is measured by the number of extracted harmonics. When the signal is fairly periodic
    %Measures, which compare the constructed cycle with the original signal, would be unhelpful due to the signal's nonharmonic and even noisy content. 

    % slope break
    The detection success of both methods is affected by signal irregularities to varying degree and in different ways. The slope break can be an indicator for signal complexity, i.e. the more modes the more complex the signal, as well as setting heuristic noise thresholds. We showed, that this is an unreliable indicator for NLSA due to its sensitivity to noise, e.g. Figs.~\ref{fig:exResolved}. 

    % noise
    SSA appears more robust to noise, as it can detect at least the fundamental oscillation even in very noisy signals where NLSA fails, Figs.~\ref{fig:exResolved}(a) and \ref{fig:exUnresolved}(a). On the other hand, noise filters do not improve the detection with SSA, Figs.~\ref{fig:exResolved}(b) and \ref{fig:exUnresolved}(b). At the same time, NLSA can fail to detect even the fundamental oscillation, Figs.~\ref{fig:exResolved}(a) and \ref{fig:exUnresolved}(a). Bandpass filter can improve detection with NLSA (five harmonics) provided the time series is regular enough, Fig.~\ref{fig:exResolved}(b), and is not affect by broadband frequency content, Fig.~\ref{fig:exUnresolved}.
    In the case of NLSA, this success can be attributed to the imposition of regularity, as shown in.\cite{giannakis2012nonlinear} This regularization process establishes connections between different euclidean  scales, approximating the full Laplacian through a diffusion kernel. Consequently, the presence of high-frequency variability poses a challenge to this enforced regularity, and influencing (and at times impeding) the detection of the overarching harmonics. This phenomenon contributes to the scarcity of instances where only the fundamental is detected. In contrast, for SSA, the major axis of the covariance ellipsoid corresponds to the fundamental—a global scale that remains unaffected by the high-frequency variability. While subsequent orders are not readily accessible, their absence does not impede the prior computation of the fundamental.

    % nonstationary behaviour
    Time series with aperiodic behavior or nonstationary dynamics involving amplitude changes can pose a challenge to our approach. We expected the detection of harmonics to fail in these cases. Instead, we find a few instances where SSA and/or NLSA extract harmonics, such as in soil water content time series (Fig.~\ref{fig:harmonics} and~\ref{fig:exAppendix3}). These seasonal cycles are difficult to interpret, as they have little in common with the original signal, Fig.~\ref{fig:exDeficient} and \ref{fig:exAppendix3}. Additional validation approaches need to be explored in such cases.

    % conclusion
    We conclude that NLSA is better suited than SSA to detect harmonic oscillation.  We propose to use NLSA's inability to detect harmonics in a signal as an indicator of the signal's deficiency to automate the selection of data that require additional processing or inspection.

    \subsection{Potential and Outlook}
    \label{sec:potentialOutlook}
    %\begin{enumerate}
        %\item why not machine learning
        %A good alternative to ML as the data is nowhere near enough for training a NN.dimension reduction is a great decomposition tool that is not as costly as e.g. such as methods that learn features in a more involved, iterated way and offers a lot of information structures that need to be accounted for and offer interpretability. 
        %\item Link Weather and climate Data: What different oscillations can be identified in dimension reduction? Specifics: subharmonics link to ENSO, nonlin DS: drift. daily cycle? 
        %\item concrete example 1/2 time series and compare seasonal cycle for statistical differences
        %\item Methodologically: How does the harmonic time-structure play out in spatio-temporal data?
        %\item outlook to other periodic variable structures: many important examples (eg cardio): versorgende medizin hin zu praeventive medizin
    %\end{enumerate}

    % land-atmosphere interactions and energy transfer
    The higher order harmonics found across sites and variables reveal a common characteristic. Thus, our data driven analysis potentially indicates energy transfer through harmonic oscillations from ecosystem drivers (air temperature and short wave radiation) to the vegetation response (land-atmosphere fluxes: GPP, NEE, $\reco$), and soil dynamics (soil water content and soil temperature).\cite{gentine2011harmonic} This suggests an intricate relationship among these ecosystem components, which could be analyzed using a multi-dimensional SSA or NLSA to reveal response delays or drift, a task reserved for future studies. Moreover, our approach could be exploited to analyze the land-atmosphere interactions during the growing season, e.g. the start and end of season.\cite{panwar2023methodological,Mora2024Macrophenology} 

    % relationship of harmonics of ground and remote sensed observation/citizen science
    Harmonic characteristics of flux variables could be employed to explore the relationship between these flux ground measurements, which have sparse spatial coverage, and observations with denser spatial coverage, such as citizen science data\cite{Mora2024Macrophenology}, or remote sensed observations, such as satellite data.\cite{martinuzzi2023learning,montero2024recurrent} 
    This could give further insights into ecosystem variability and impact of climate change on vegetation\cite{ghil2020physics} as well as vegetation-climate feedback loops.\cite{mahecha2024biodiversity}
    %, or land-atmosphere model data.

    Our approach is widely applicable. For example, Higher order harmonics have also been detected with these methods in climate observations, e.g. four harmonics were detected in intra-seasonal characteristics in satellite infrared brightness temperature\cite{szekely2016extraction} with NLSA, and seven harmonics were detected in Sentinel-1 Interferometric Synthetic Aperture Radar data with multichannel SSA\cite{walwer2023multichannel}, which is a spatio-temporal SSA approach. 

    % subharmonics
    A spatio-temporal Ansatz with these methods could explore the link between subharmonics of land-atmosphere carbon fluxes, oscillations with periods of multiple years, and global climate patterns such as El Ni\~{n}o-Southern Oscillation\cite{dijkstra2013nonlinear} (ENSO). Our analysis only tentatively suggests existence of subharmonics, e.g. Fig.~\ref{fig:exResolved}. The relationship between the seasonality in sea surface temperature and ENSO have been studied\cite{wang2018antarctic} with NLSA. 

    % extremes and prediction
    Extreme heatwaves and droughts can affect carbon gas exchange, indicating reduced vegetation productivity, as observed in GPP during the summer of 2018. These impacts seem to be distinguishable from the harmonic oscillations of the seasonal cycle. We conjecture that nonharmonic modes could represent such extreme events. \cite{giannakis2013nonlinear} Unlike SSA, NLSA is shown to successfully detect such intermittent modes related to rare events. Ideally, we would like to predict the impact of extreme events on vegetation. However, prediction is challenging for such methods and remain a challenge even for sophisticated machine learning models.\cite{martinuzzi2023learning}
    
%#########################
%   S E C T I O N : Conclusion
%########################
\section{Conclusion}
\label{sec:conlusion}
Land-atmosphere exchange of CO$_2$ is complex. The underlying processes are characterized by dynamics which are not necessarily dominated by large variance. 
To extract the corresponding seasonal dynamics we compared two data driven methods: NLSA a graph-based algorithm, and SSA which identifies dynamic characteristics by maximizing the data's variance.
We showed that NLSA is better suited than SSA to extract higher order harmonics, which are necessary to represent the seasonal cycle accurately.
When extraction of harmonics was not possible, NLSA proved to be a more reliable indicator of nonstationary dynamics than SSA.
%, even if it is more robust to artifacts. 

% OLD TEXT
%\subsection{general conclusion}
%we extract complicated dynamics of unharmonic nonlinear time series variables reliably with a high resolution.

%\subsection{opinion: climate science perspective}
%highly recommend to use NLSA instead of SSA: it ties the harmonic structure to the fundamental frequency: this yields reliability in terms of not resolving an too-irregular signal both due to high frequency variability or seasonal cycle disturbances. Very little is gained in climate science by robustly accessing the fundamental.
%
% NLSA is great, here is why.
% filter noise for best results.
% We propose to employ NLSA's feature of failing ...
% refer to different algorithm principles, i.e. SSA always works as it maximises variance, whereas NLSA works because .
% (along the lines of ill posed problem).
% No results due to lack of periodic dynamics, SSA almost always yields results due to the computation principle.
%link to figures

\begin{acknowledgments}

We thank Holger Kantz and Markus Reichstein for insightful discussions. This work was initiated during a workshop for the proposed excellence cluster Breathing Nature. We thank the Breathing Nature community for the interdisciplinary exchange. KM acknowledges funding by the Saxon State Ministry for Science, Culture and Tourism (SMWK) – [3-7304/35/6-2021/48880]. We thank the German Federal Ministry for Economic Affairs and Climate Action for supporting us via the ML4Earth project (grant number: 50EE2201B).
\end{acknowledgments}

\section*{Author Declarations}

\subsection*{Conflict of Interest}

The authors have no conflicts to disclose.

\subsection*{Author Contributions}

\textbf{Karin Mora:} 
Conceptualization (lead); 
Methodology (lead);
Formal analysis (equal);
Project administration (lead);
Supervision (lead);
Writing - original draft (lead);
Writing - review \& editing (equal).

\textbf{Leonard Schulz:} 
Conceptualization (supporting);
Methodology (supporting);
Formal analysis (equal);
Software (lead); 
Writing - original draft (supporting);
Writing - review \& editing (supporting).

\textbf{Jürgen Vollmer:} 
Conceptualization (supporting);
Methodology (supporting);
Formal analysis (supporting);
Project administration (supporting);
Supervision (supporting);
%Writing - original draft (supporting);
Writing - review \& editing (equal).

\textbf{Miguel D. Mahecha:} 
Conceptualization (supporting);
Methodology (supporting);
Formal analysis (supporting);
Supervision (supporting);
%Writing - original draft (supporting);
Writing - review \& editing (equal).

\section*{Data Availability}
This study used the ICOS 2020 Warm Winter dataset\cite{WarmWinter2020Team2022} 
(\href{https://doi.org/10.18160/2G60-ZHAK}{https://doi.org/10.18160/2G60-ZHAK}).
The code to reproduce our results and figures is freely available.\cite{dimensionality} (\href{https://github.com/SidxA/dimensionality}{https://github.com/SidxA/dimensionality}).

%
%The data that support the findings of this study are openly available via ICOS (\url{www.icos-cp.eu}) for Europe or FLUXNET at the global scale (\url{fluxnet.org}).

\bibliography{seasonalCycle}
\appendix

% \section{Appendixes}

\section{Data analysis workflow details}
\label{app:methodsAll}

    \subsection{Dimension reduction}
    \label{app:dimred}
    Dimension reduction decomposes data into adaptive, orthogonal patterns (modes).
    In the case of univariate observations with a time ordering $x_t$,
    modes of specified length $W$ are obtained.
    Creating information redundancy in a high dimensional embedding space $X_t$,
    a finite set of $k$ projection directions $X_t^k = \sum_1^k u_i \sigma_i v_i(t)$ is chosen to optimize specific conditions.
    The modes $u_i$ and the projections onto them $\sigma_i v_i(t)$ correspond to data-driven orthogonal basis functions.
    
    \subsection{Uniform delay embedding}
    \label{app:embedding}
    From the vector $x_t \in \mathbb R^N$ containing the time series of $N$ uniformly sampled measurements,
    a delay trajectory matrix $X \in \mathbb R^{W \times P}$ is constructed by uniform delay embedding,
    
    \begin{equation}
    x_t \mapsto X = 
    \begin{bmatrix}
    x_{1}&, \hdots &,x_{W} \\
    x_{2}&, \hdots &,x_{W+1} \\
    &\vdots& \\
    x_{N-W}&, \hdots &,x_{N}
    \end{bmatrix}.
    \label{eq:embedding}
    \end{equation}
    
    $W \in \{1,\hdots,N/2\}$ is the delay window length and $P = N - W + 1$ the number of thus created windows.
    This matrix corresponds to a manifold of points (delay windows) or hypersurface in $\mathbb R^W$.
    Here, the longest full-year delay embedding time of $W = \SI{7}{a}$ is used as parameter choice to guarantee the best possible spectral resolution of the modes with $N = \SI{14}{a}$.
    
    \subsection{Singular Spectrum Analysis}
    \label{app:ssa}
    Based on the Singular Value decomposition
    
    \begin{equation}
    X = U \Sigma V^T,
    \end{equation}
    
    Singular Spectrum Analysis\cite{broomhead1986extracting} (SSA) utilizes the $k$ first singular values and left singular vectors of $X$.
    These Empirical Orthogonal Functions\cite{bjornsson1997manual} (EOF) correspond to the axes of the linear covariance ellipsoid $XX^T$ and therefore best decompose $X$ in terms of squared reconstruction error $
    \mathbb E_2 = \sum_t | X_t - \sum_i \sigma_i u_i |^2$.
     
    \subsection{Nonlinear Laplacian Spectral Analysis}
    \label{app:nlsa}
    
    SSA does not pose any constraints on the shape of the modes except for the orthogonal structure.
    Nonlinear Laplacian Spectral Analysis\cite{giannakis2012nonlinear} (NLSA) additionally enforces regularity on the projections to recover local nonlinear behavior.
    The first $k$ eigenfunctions $\phi_i$ of the discrete Laplace-Beltrami operator $\Delta$ on the manifold are approximated.
    This best reconstructs the data in terms of the heat kernel $\exp(-t\Delta)$, which is equivalent to the Riemannian metric, i.e. the inner product of the manifolds tangent vectors.\cite{giannakis2012nonlinear}
    It is shown\cite{coifman2006diffusion} that the $\phi_i$ can be sufficiently approximated by the eigenfunctions of a renormalized diffusion kernel.
    The diffusion kernel 
    
    \begin{equation}
    J_\epsilon(z) = \exp(-z^2/2\epsilon),\,\, z= |x-y| \text{ at scale } \sqrt\epsilon
    \end{equation}
    
    stores the diffusion distances between all manifold points.
    This is helpful to approximate the sampling density $\rho$ on the manifold
    
    \begin{equation}
    \mathsf Z(x) =  \sum_{i=1}^P J(x,x_i) \approx \int_X J(x,y)  \rho(y)  dy.
    \end{equation}
    
    A discrete representation of the Laplace-Beltrami operator on the manifold is achieved by this anisotropic renormalization\cite{jost2008riemannian},
    
    \begin{equation}
    T = \frac{\mathsf{Z}^{-1} J \mathsf{Z}^{-1}}{\text{diag}(\mathsf{Z}^{-1} J \mathsf{Z}^{-1})}.
    \label{eq:anisotropic}
    \end{equation}
    
    The projections $X u_i = \sigma_i v_i$ are regular,
    because the $u_i$ are directly linked to the Riemannian metric.
    The approximated eigenfunctions of $T$ correspond to the EOF.
    The kernel scale parameter $\sqrt \epsilon$ is crucial for the performance of the algorithm.
    Its size is linked to the assumed sampling density and determines the connectivity of the manifold graph: for smallest $\epsilon$, no vertices and for largest $\epsilon$, all nodes are connected. Either case does not yield any information about the local structure.
    We suggest an simple automated choice,
    performed by a sampling approach on the embedding manifold:
    \begin{enumerate}
    \item for multiple values of $\epsilon$ compute the weight $\mathsf Z^c = \sum_{i,j}^c J(x_i,x_j)$ for a random manifold subset indexed by $1,\hdots,c$
    \item estimate the turning point $\mathsf Z^c_t$ of the resulting bidirectional kernel saturation curve using a tanh fit by least square rgression
    \item choose $\epsilon$ such that $\mathsf Z^c \approx \mathsf Z^c_t/e$ 
    \item take the median of multiple sampling runs for $\epsilon$
    \end{enumerate}
    
    The other parameters $\alpha = 1,t=1$ of the \textsl{Diffusion Map} algorithm for general diffusion kernel computations are not considered here.
    
    \subsection{Time series reconstruction}
    \label{app:TSrecon}
    
    The Uniform Delay Embedding represents the time steps at the beginning and end of the original time series less often in the embedding space, due to the limited number of covering windows.
    Reconstructing the time RC series from the recovered reduced dimensions $u_i$ and the corresponding projections PC takes this trapezoid shape into account.\cite{ghil2002advanced}
    
    \begin{eqnarray}
    \text{RC}_\kappa(t) = \mathtt M(t)^{-1} \sum_i^\kappa \sum_{j=\mathtt L(t)}^{\mathtt U(t)} PC_i(t-j+1) u_i(j)\text{, with} \nonumber\\
    (\mathtt M,\mathtt L, \mathtt U)(t) =
    \begin{cases}
    (t,1,t) & 1 \leq t \leq W-1 \\
    (W,1,W) & W \leq t \leq P \\
    (N-t+1,t-N+W,W) & P+1 \leq t \leq N
    \end{cases}
    \end{eqnarray}
    
    \subsection{Signal centralization}
    \label{app:centralization}
    
    The covariance-relation of the SSA requires the data to be centered at zero mean.
    This is performed on the embedding space.
    
    \begin{equation}
    \tilde X_i = (X_i - \hat X_i) / \sigma^2(X_i)
    \end{equation}
    
    Also, centralizing the time series directly allows to link the individual projection scales (diagonal entries of $\Lambda$) from the projections $u_i \in \text{EOF}$ directly to corresponding variance of the time series component,
    
    \begin{equation}
    \Lambda = \frac{1}{W} \text{EOF}^T \cdot X^T \cdot X \cdot \text{EOF}.
    \label{eq:var}
    \end{equation}
    
    \subsection{Orthogonality of dimensions}
    \label{app:orthogonality}
    
    Both dimensionality reductions methods, SSA and NLSA, can be expressed as eigendecompositions of respective kernel matrices.
    For SSA, the covariance kernel $X X^T$ and the renormalized Diffusion kernel $T$ for NLSA.
    The reduced dimensions are the obtained eigenvectors, a non-unique orthonormal basis that diagonalizes the kernel matrix.
    Two identical eigenvalues lead to a linearly independent and degenerate pair of eigenvectors - resolving each oscillatory component by two modes.
    Because the kernel matrices are positive and symmetric, eigenvectors corresponding to different eigenvalues are mutually orthogonal.
    In oscillatory time-series, this orthogonality creates the harmonic structure.
    Trying to resolve such signal,
    the fundamental driving frequency is the first approximation and holds the most variability.
    All dimensions that can exclusively be attributed to the periodic, non-harmonic signal need to be combined to recover what can be identified from this signal by the dimensionality reduction.
    
    \subsection{Harmonic structure detection}
    \label{app:harmonic}
    
    To identify purely harmonic modes, we propose a simple filter based on the spectral content using Fast-Fourier-Transform\cite{cooley1965algorithm} (FFT),
    \begin{equation}
    X_n = \sum_{n=0}^{N-1} x_n \exp(\frac{-2\pi}{N} k_n).
    \end{equation}
    FFT provides approximated harmonic power $X_n$ for frequencies $k_n$ 
    %in the range $\{0, \hdots,\frac{1}{N-1} \}$ (FIND THIS OUT).
    First, the information content of a mode in the frequency range of the harmonics $k = \{ 0,\hdots,7 \}$ of the fundamental frequency is fitted by a Gaussian peak,
    \begin{equation}
    g(k) = \frac1{\sigma \sqrt{2\pi}} \exp(-\frac1 2 \frac{(k-\mu_k)^2}{\sigma^2}),
    \end{equation}
    by least square regression, yielding the parameters $\sigma$ and $\mu_k$.
    If this peak resides at a harmonic frequency $\text{abs}(\mu_k-f) \leq \epsilon_f \text{ for any } f \in \{ 1,2,3,4 \} \unit{\year^{-1}}$ and the largest remaining spectral power is not larger then a portion $\epsilon_p$ of this peak $\frac1{\sigma \sqrt{2\pi}} $, the dimension is classified as harmonic.
    Both parameters are chosen to be $\epsilon_f = \epsilon_p = 0.15$.
    
    We apply a low-pass filter to the time series to analyze the presence of high frequency variability on the dimension reduction performance.
    We use the Inverse-Fast-Fourier-Transform (IFFT),
    \begin{equation}
    x_n = \frac 1 N \sum_{n=0}^{N-1} X_n \exp(\frac{i2\pi}{N} k_n)
    \end{equation}
    setting $X_k = 0 \quad\forall\quad k_n \geq f_l$, therefore omitting all spectral content above $f_l = \{3,4,6\} \unit{\year^{-1}}$ respectively.
    
    \subsection{Data characterization measures}
    \label{app:dataChar}
    We compute the relative \textbf{extent of high frequency variability}, that can impede the extraction of the seasonal cycle from the time series.
    The sum of spectral power $\sum X(k) \quad\forall\quad k \geq \SI{6}{\year^{-1}}$, normalised by $\sum_k X_k$.
    This frequency threshold provides a good compromise of information loss in this data(harmonics up to that point are resolved).
    % Code: 
    
    The relative power by FFT of the harmonic components belonging to the fundamental annual frequency $\sum X_n \text{ where } k_n = \{1,2,3,4\} \unit{\year^{-1}}$, normalized by$\sum_k X_k$ yields a measure for the \textbf{time series regularity}.
    % Code: 05_processing line 331
      
    In information theory the \textbf{sample entropy}\cite{richman2000physiological} (SampEn) estimates the complexity by approximating the randomness in a time series.
    %, referring to the non-Gaussianity of time series.
    Here we choose the embedding dimension $m=2$, and the margin of tolerance $r = 0.2 \text{std}$, where std is the standard deviation of the whole time series.
    %and Euclidean norm.
    The sample entropy $SampEn$ is defined as
    \begin{equation}
    SampEn = - \ln \frac A B,
    \end{equation}
    which yields the negative logarithm of the proportion of $A$ -- the number of subsequent delay-3-windows being not farther apart then $0.2\text{std}$, and $B$ -- the number of subsequent delay-2-windows being no farther apart then $0.2\text{std}$.

    To display the data metrics, a \textbf{binning scheme} (histogram) with 3 bins is employed across all time series.
    The difference between the smallest and largest value is divided into 3 equally spaced intervals.
    The relative power of high frequency variability is only binned across unfiltered time series.
    
    The quality of measurements and gap-filling procedures is indicated by additional quality flag variables.\cite{baldocchi2001fluxnet,baldocchi2018inter}
    To detect longer persisting issues, we classify \textbf{persistent QF} as the existence of windows of length $\geq \frac 1 2 \text{a}$, where all quality flags are below whole-series-mean.
    
    \subsection{Computational cost}
    \label{app:compCost}
    
    The total computational cost depends on the length $N$ of the time series, the delay embedding parameter $W$ and the corresponding number of delay embedding data points $P$.\cite{cook1983overview}
    Here, $N=5114,W=2556,P=2559$.
    Apart from generating the delay embedding matrix,
    the SVD to calculate SSA is $O(W P \text{min}(W, P))$.
    The covariance approach is less expensive with $O(W^2)$.
    In NLSA, the distance evaluation is costliest with $O(P^2)$,
    which can be accelerated by using graph-theoretical K-n-nearest neighbor tree distances to $O(P \text{log} P)$.
    Sampling iterations to approximate $\epsilon$ cost $O(P)$.
    
%\newpage
\section{Supplementary Figures}
%\label{chap:AppendixData}
\label{app:suppFig}

\newpage
\subsection{Harmonic extraction across all sites and variables: number of harmonics}
\label{sec:harmonicsNum}
\begin{figure}[h!]
\centering
\includegraphics[width=0.9\linewidth]{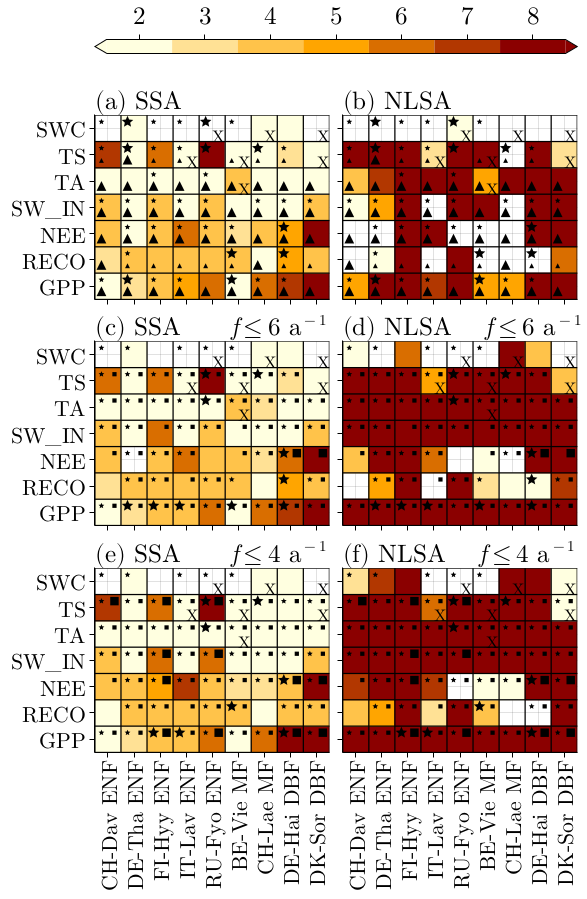}
\caption{
        Harmonic oscillation extraction across variables and sites: number of harmonic modes (color bar) detected in Fig.~\ref{fig:harmonics}.
        Panels: SSA (a,c,e), NLSA (b,d,f).
        Sites (Table~\ref{table:sites}) with land cover types ($x$-axis): evergreen needleleaf forest (ENF), mixed forest (MF), deciduous broadleaf forest (DBF).
        Variables ($y$-axis): soil water content (SWC), soil temperature (TS), air temperature (TA), Short wave radiation (SW\_IN), net ecosystem exchange (NEE), ecosystem respiration ($\reco$), gross primary productivity (GPP).
        Filters: 
        unfiltered (a,b), 
        low-pass filter with $f = \SI{6}{\year^{-1}}$ (c,d) 
        and with $\SI{4}{\year^{-1}}$ (e,f). 
        Markers for time series properties:
        X = quality flags of at least $\text{a}/2$ length,
        $\blacktriangle$ = high frequency variability (high, low, none),
        $\star$ = sample entropy (high, low, no randomness),
        $\blacksquare$ = regularity (high, low, none). 
        }
\label{fig:harmonicsNum}
\end{figure}

%\newpage
\subsection{Additional examples for Cases~1--4}
\begin{figure*}[htb!]
    \centering
    \includegraphics[width=\linewidth]{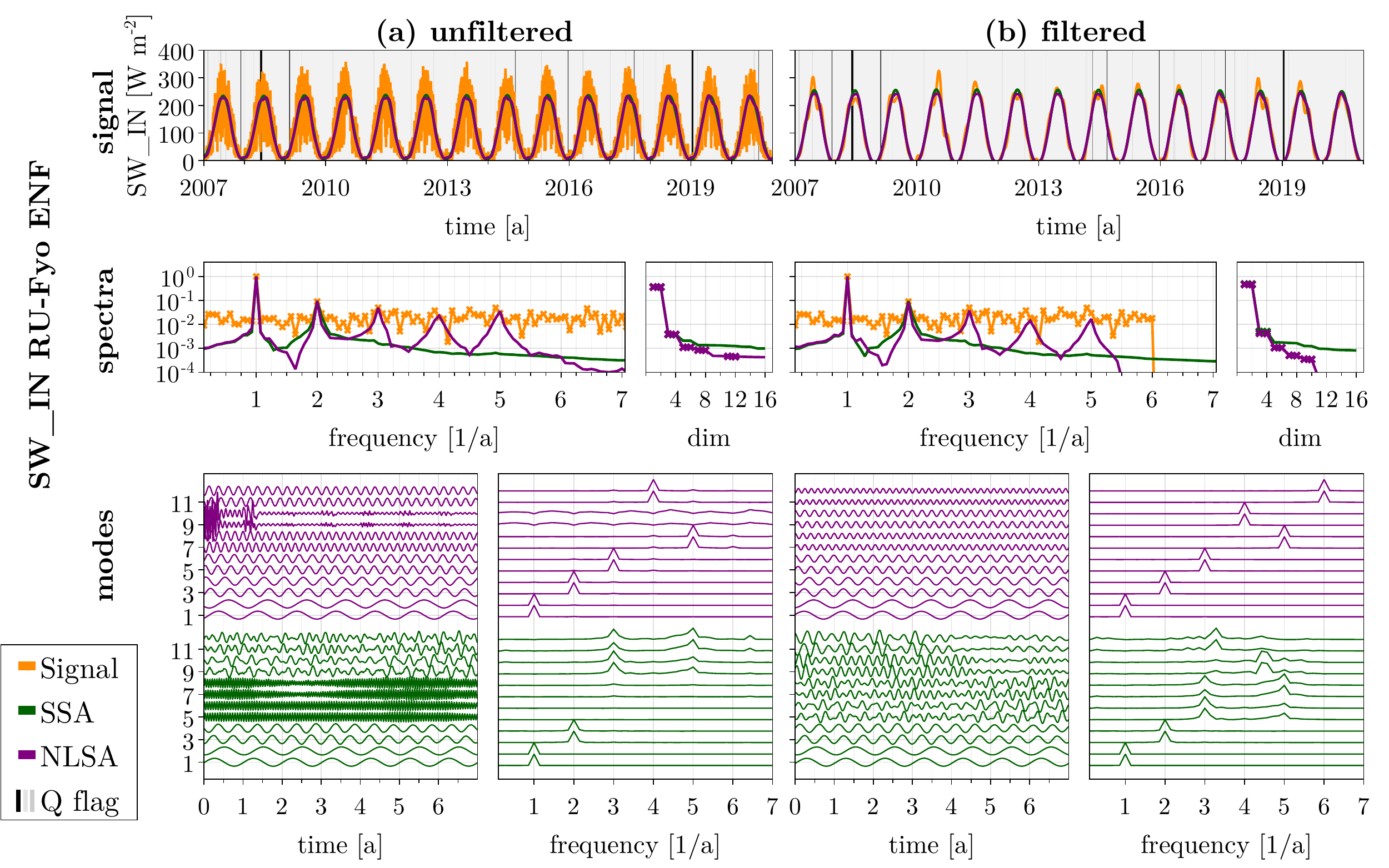}
    \caption{
    Analysis of a regular time series with isolated quality flags: incoming shortwave radiation power (SW\_IN) of an evergreen needleleaf forest (ENF) in Western Russia (Fyodorovskoye site, RU-Fyo ENF). 
        \\
        Panels: Analysis of an unfiltered (a) and filtered (b) measurement signal. Each panel column illustrates the following analysis parts.
        Panel \emph{signal}: signal and the corresponding constructed seasonal cycle with SSA or NLSA, if harmonics are detected. The quality flag (QF) is indicated by a gray scale, i.e. lowest signal quality in black.
        Panel \emph{spectra}: relative FFT power spectrum (left) of the signal and the SSA and NLSA constructed seasonal cycle, and the dimension reduction spectrum (right).
        Panel \emph{modes}: shapes of the 12 most dominant modes (left) and their corresponding spectra (right). \\
        Main result: Both SSA and NLSA extract higher order harmonics, NLSA up to fifth, and more after filtering. SSA modes exhibit information overlap between third and fifth harmonic order.
    }
    \label{fig:exAppendix1}
\end{figure*}

\begin{figure*}[htb!]
    \includegraphics[width=\linewidth]{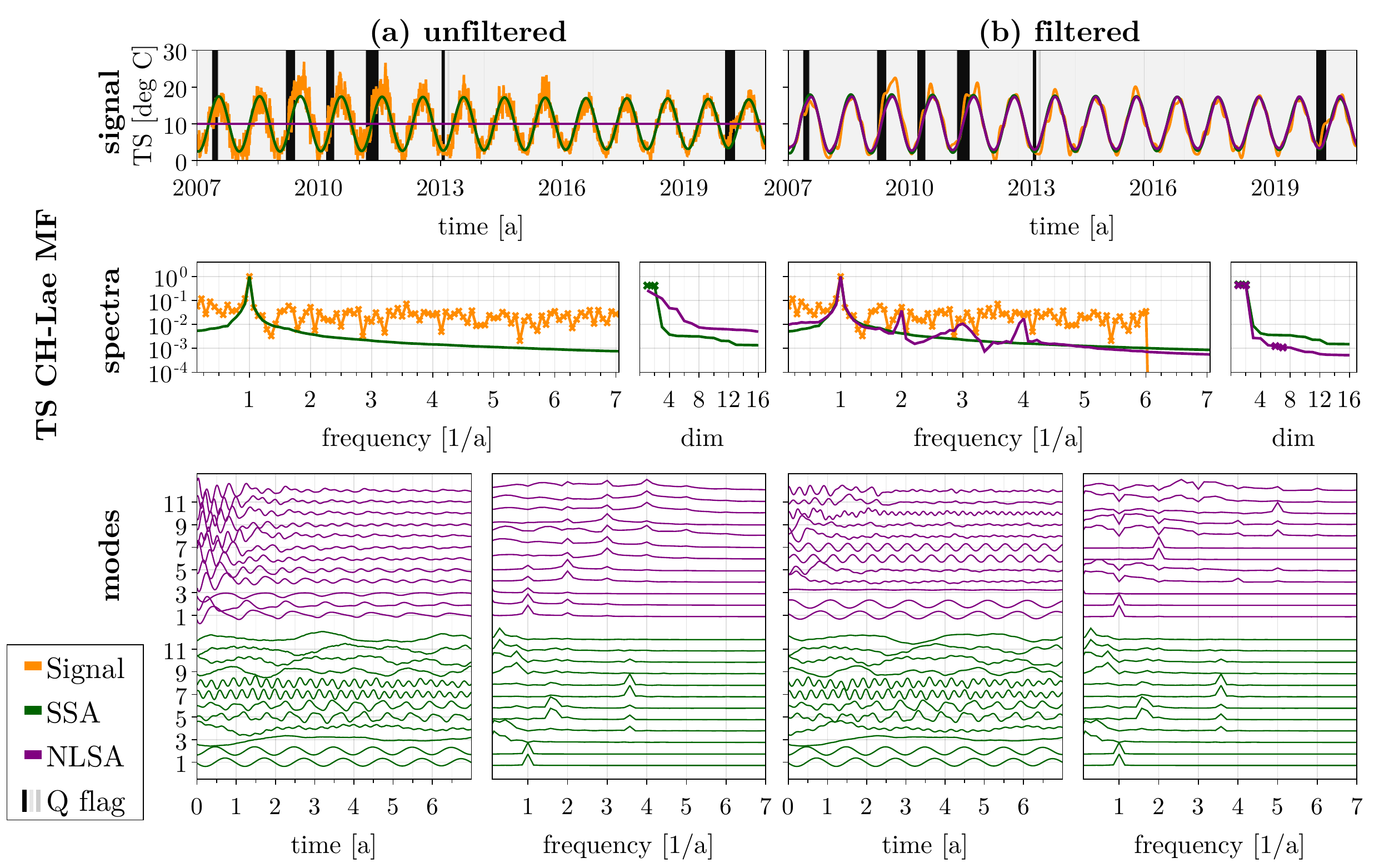}
    \caption{
    %(a) SSA computed fundamental. NLSA did not compute harmonic structure. Instead, decreasing amplitude modulation within modes hinting at a information resolution conflict between the seasonal cycle and the changes in high frequency variability. Modes show distributed harmonic orders. 
    %(b) (b) Time series shows summer amplitude deviations in 2009,2010,2011 from otherwise regular appearance. NLSA yields four harmonic orders.
    %%
    Analysis of a regular time series with high frequency variability and isolated quality flags: Soil temperature (TS) of a mixed forest (MF) in Northern Switzerland (Laegeren mountain site, CH-Lae MF).
        \\
        Panels: Analysis of an unfiltered (a) and filtered (b) measurement signal. Each panel column illustrates the following analysis parts.
        Panel \emph{signal}: signal and the corresponding constructed seasonal cycle with SSA or NLSA, if harmonics are detected. The quality flag (QF) is indicated by a gray scale, i.e. lowest signal quality in black.
        Panel \emph{spectra}: relative FFT power spectrum (left) of the signal and the SSA and NLSA constructed seasonal cycle, and the dimension reduction spectrum (right).
        Panel \emph{modes}: shapes of the 12 most dominant modes (left) and their corresponding spectra (right). \\
        Main result: SSA extracts the fundamental oscillation, while NLSA extracts forth order harmonics after filtering. Before filtering high frequency variability affects signal separation for NLSA.
    }
    \label{fig:exAppendix2}
\end{figure*}

\begin{figure*}[htb!]
    \includegraphics[width=\linewidth]{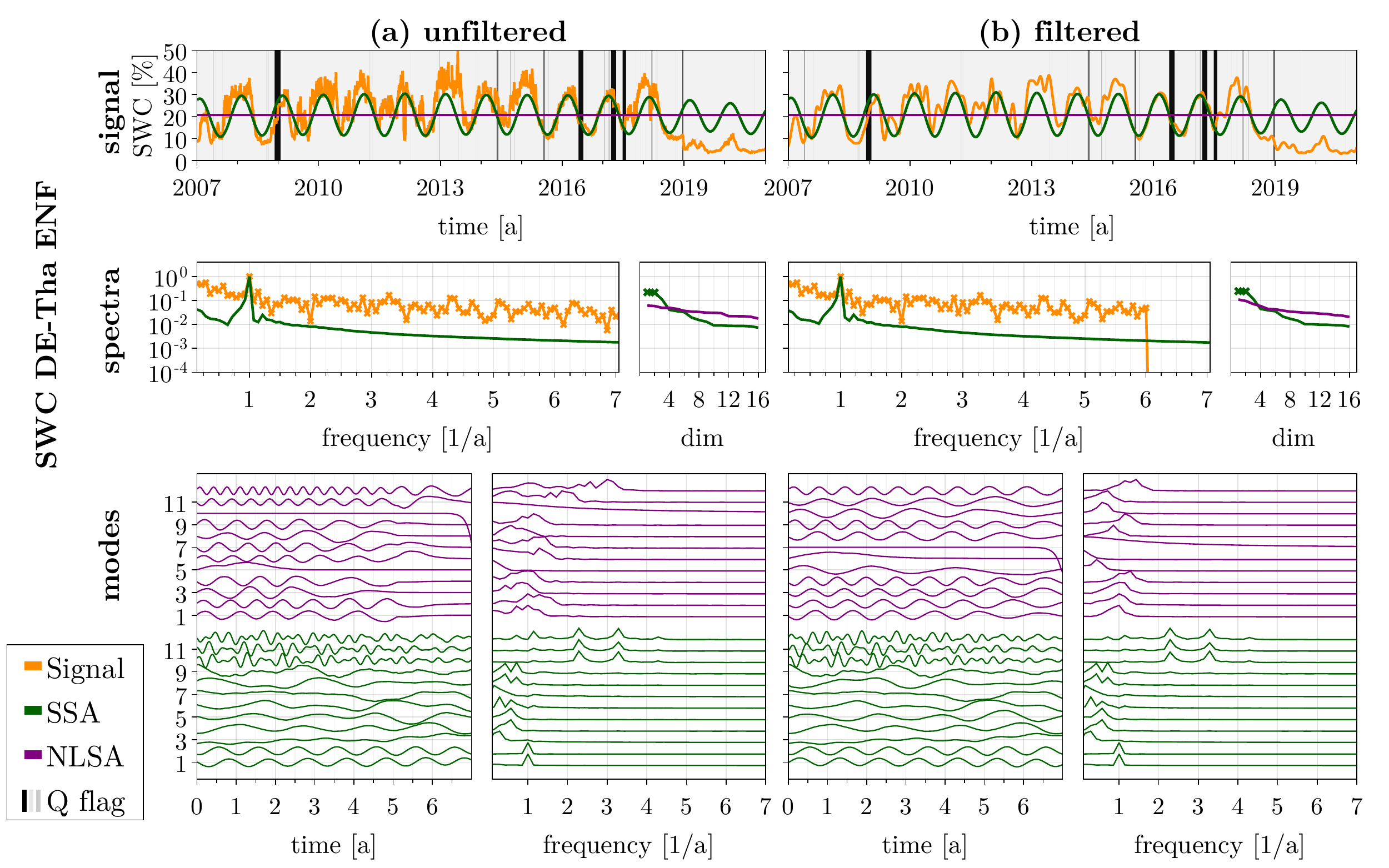}
    \caption{
    %Irregular soil water content time series with high frequency variability and abrupt changes in the seasonal cycle from 2018 on. Identical layout as Fig.~ \ref{fig:exRegular}. (a) SSA identified the fundamental solely. NLSA did not compute harmonic structure, but frequency modulated modes clustered around it.
    %(b) Times series appear irregular. Minor changes in mode spectral composition for both methods, towards harmonic structure for NLSA (i.e. less frequency modulation) - not detected still.
    Analysis of an irregular time series with amplitude change and isolated quality flags: soil water content (SWC) of an evergreen needleleaf forest (ENF) in Eastern Germany (Tharandt forest site, DE-Tha ENF).
        \\
        Panels: Analysis of an unfiltered (a) and filtered (b) measurement signal. Each panel column illustrates the following analysis parts.
        Panel \emph{signal}: signal and the corresponding constructed seasonal cycle with SSA or NLSA, if harmonics are detected. The quality flag (QF) is indicated by a gray scale, i.e. lowest signal quality in black.
        Panel \emph{spectra}: relative FFT power spectrum (left) of the signal and the SSA and NLSA constructed seasonal cycle, and the dimension reduction spectrum (right).
        Panel \emph{modes}: shapes of the 12 most dominant modes (left) and their corresponding spectra (right). \\
        Main result: SSA extracts the fundamental oscillation, which does not align with the time series from 2018. NLSA does not extract pure oscillations indicating the signal's nonstationary dynamics.
    }
    \label{fig:exAppendix3}
\end{figure*}

\begin{enumerate}
    \item[Fig.~\ref{fig:exAppendix1}] additional example for Case~1.
    %: Regular time series. 
    \item[Fig.~\ref{fig:exAppendix2}] additional example for Case~2.
    %: Time series with high frequency variability. 
    \item[Fig.~\ref{fig:exAppendix3}] additional example for Cases~3 and 4.
    %: Time series with broadband variability and amplitude change.
\end{enumerate}

%\newpage

%\newpage
%\subsection{Additional example for Case~2: Time series with high frequency variability}
\vfill

%\newpage
%\subsection{Additional example for Cases~3 and 4: Time series with broadband variability and amplitude change}

\end{document}